\documentclass{article}
\usepackage{geometry}
\usepackage[utf8]{inputenc} % allow utf-8 input
\usepackage[T1]{fontenc}    % use 8-bit T1 fonts
\usepackage{hyperref}       % hyperlinks
\usepackage{url}            % simple URL typesetting
\usepackage{booktabs}       % professional-quality tables
\usepackage{amsfonts}       % blackboard math symbols
\usepackage{nicefrac}       % compact symbols for 1/2, etc.
\usepackage{microtype}      % microtypography

\usepackage[square, compress]{natbib}

\usepackage[ruled, vlined, linesnumbered]{algorithm2e}

\usepackage{amsfonts}
\usepackage{amsthm,amsmath}
\usepackage{amssymb}
\usepackage{bbm}
\usepackage{bm}
\newtheorem{lemma}{Lemma}

\newtheorem{proposition}{Proposition}
\usepackage{dsfont}

\usepackage{hyperref}
\usepackage{enumitem}
\usepackage{comment}

\usepackage{adjustbox}
\usepackage{microtype}
\usepackage{graphicx}
\usepackage{subfigure}
\usepackage{booktabs}

\usepackage{color}
\usepackage[dvipsnames]{xcolor}

\definecolor{cgreen}{HTML}{006600}

\usepackage{xspace}
\usepackage{booktabs}
\usepackage{multirow}
\newcommand{\ra}[1]{\renewcommand{\arraystretch}{#1}}
\newtheorem{Definition}{Definition}
\newtheorem{Example}{Example}

\newcommand\IR{\mathbb{R}}
\newcommand\IE{\mathbb{E}}
\newcommand\IP{\mathbb{P}}

\newcommand\calO{\mathcal{O}}
\newcommand\calI{\mathcal{I}}
\newcommand\OUI{\mathcal{O} \cup \mathcal{I}}
\newcommand\calS{\mathcal{S}}

\usepackage{graphicx}
\usepackage{lscape}
\usepackage{balance}

\title{Robust Kernel Density Estimation \\ with Median-of-Means principle}

\author{%
  Pierre~Humbert\thanks{These authors contributed equally to this work.} \quad
 Batiste~Le Bars$^*$ \quad
 Ludovic~Minvielle$^*$ \vspace{.5em} \\ \vspace{.5em}
 Nicolas~Vayatis \\
Université Paris-Saclay, CNRS, ENS Paris-Saclay, Centre Borelli \\ F-91190 Gif-sur-Yvette, France
}
\date{June 5, 2020}

\begin{document}
\maketitle
\begin{abstract}
In this paper, we introduce a robust nonparametric density estimator combining the popular Kernel Density Estimation method and the Median-of-Means principle (MoM-KDE). This estimator is shown to achieve robustness to any kind of anomalous data, even in the case of adversarial contamination. In particular, while previous works only prove consistency results under known contamination model, this work provides finite-sample high-probability error-bounds without \textit{a priori} knowledge on the outliers. Finally, when compared with other robust kernel estimators, we show that MoM-KDE achieves competitive results while having significant lower computational complexity.
\end{abstract}

\section{Introduction}

Over the past years, the task of learning in the presence of outliers has become an increasingly important objective in both statistics and machine learning. Indeed, in many situations, training data can be contaminated by undesired samples, which may badly affect the resulting learning task, especially in adversarial settings.
Building robust estimators and algorithms that are resilient to outliers is therefore becoming crucial in many learning procedures. In particular, the inference of a probability density function from a contaminated random sample is of major concerns.

Density estimation methods are mostly divided into parametric and nonparametric techniques.
Among the nonparametric family, the Kernel Density Estimator (KDE) is probably the most known and used for both univariate and multivariate densities \citep{parzen1962estimation, silverman1986density, scott2015multivariate}, but it also known to be sensitive to dataset contaminated by outliers \citep{kim2011robustness, kim2012robust, vandermeulen2014robust}.
The construction of robust KDE is therefore an important area of research,
that can have useful applications such as anomaly detection and resilience to adversarial data corruption.
Yet, only few works have proposed such robust estimators.

\citet{kim2012robust} proposed to combine KDE with ideas from M-estimation to construct the so-called Robust Kernel Density Estimator (RKDE). However, no consistency results were provided and robustness was rather shown experimentally. Later, RKDE was proven to converge to the true density, however at the condition that the dataset remains uncorrupted \citep{vandermeulen2013consistency}.
More recently, \citet{vandermeulen2014robust} proposed another robust estimator, called Scaled and Projected KDE (SPKDE).
Authors proved the $L_1$-consistency of SPKDE under a variant of the Huber's $\varepsilon$-contamination model where two strong assumptions are made \citep{huber1992robust}. First, the contamination parameter $\varepsilon$ is known, and second, the outliers are drawn from an uniform distribution when outside the support of the true density.
Unfortunately, as they did not provided rates of convergence, it still remains unclear at which speed SPKDE converges to the true density.
Finally, both RKDE and SPKDE require iterative algorithms to compute their estimators, increasing the overall complexity of their construction.

In statistical analysis, another idea to construct robust estimators is to use the Median-of-Means principle (MoM).
Introduced by \citet{nemirovsky1983problem}, \citet{jerrum1986random}, and \citet{alon1999space}, the MoM was first designed to estimate the mean of a real random variable. It relies on the simple idea that rather than taking the average of all the observations, the sample is split in several non-overlapping blocks over which the mean is computed. The MoM estimator is then defined as the median of these means.
Easy to compute, the MoM properties have been studied by \citet{minsker2015geometric} and \citet{devroye2016sub} to estimate the means of heavy-tailed distributions.
Furthermore, due to its robustness to outliers, MoM-based estimators have recently gained a renewed of interest in the machine learning community  \citep{lecue2018robust,lecue2019learning}.

\paragraph{Contributions.} In this paper, we propose a new robust nonparametric density estimator based on the combination of the Kernel Density Estimation method and the Median-of-Means principle (MoM-KDE). We place ourselves in a more general framework than the classical Huber contamination model, called $\OUI$, which gets rid of any assumption on the outliers. We demonstrate the statistical performance of the estimator through finite-sample high-confidence error bounds in the $L_\infty$-norm and show that MoM-KDE's convergence rate is the same as KDE without outliers. Additionally, we prove the consistency in the $L_1$-norm, which is known to reflect the global performance of the estimate. To the best of our knowledge, this is the first work that presents such results in the context of robust kernel density estimation, especially under the $\OUI$ framework.
Finally, we demonstrate the empirical performance of MoM-KDE on both synthetic and real data and show the practical interest of such estimator as it has a lower complexity than the baseline RKDE and SPKDE.

\section{Median-of-Means Kernel Density Estimation}
We first recall the classical kernel density estimator. Let $X_1, \cdots, X_n$ be independent and identically distributed (i.i.d.) random variables
that have a probability density function (pdf) $f(\cdot)$ with respect to the Lebesgue measure on $\IR^d$. The Kernel Density Estimate of $f$ (KDE), also called the \textit{Parzen–Rosenblatt estimator}, is a nonparametric estimator given by
\begin{equation}
    \label{eq:kde}
    \hat{f}_n(x) = \dfrac{1}{nh^d} \sum^n_{i=1} K\left(\dfrac{X_i - x}{h}\right) \;,
\end{equation}
where $h>0$ and $K : \IR^d \longrightarrow \IR_{+}$ is an integrable function satisfying $\displaystyle\int K(u)du = 1$ \citep{tsybakov2008introduction}. Such a function $K(\cdot)$ is called a \textit{kernel} and the parameter $h$ is called the \textit{bandwidth} of the estimator.
The bandwidth is a smoothing parameter that controls the bias-variance tradeoff of $\hat{f}_n(\cdot)$ with respect to the input data. \\

\noindent
While this estimator is central in statistic, a major drawback is its weakness against outliers  \citep{kim2008robust, kim2011robustness, kim2012robust, vandermeulen2014robust}. Indeed, as it assigns uniform weights $1/n$ to every points regardless of whether $X_i$ is an outlier or not, inliers and outliers contribute equally in the construction of the KDE, which results in undesired ``bumps'' over outlier locations in the final estimated density (see Figure \ref{fig:bump_kde}).
In the following, we propose a KDE-based density estimator robust to the presence of outliers in the sample set. These outliers are considered in a general framework described in the next section.
\subsection{Outlier setup}

Throughout the paper, we consider the $\OUI$ framework introduced by \citet{lecue2019learning}. This very general framework allows the presence of outliers in the dataset and relax the standard i.i.d. assumption on each observation. We therefore assume that the $n$ random variables are partitioned into two (unknown) groups: a subset $\{X_i \mid i \in \mathcal{I}\}$ made of inliers, and another subset $\{X_i \mid i \in \mathcal{O}\}$ made of outliers such that $\mathcal{O} \cap \mathcal{I} = \emptyset$ and $\mathcal{O} \cup \mathcal{I} = \{1,\ldots,n\}$. While we suppose the $X_{i \in \mathcal{I}}$ are i.i.d. from a distribution that admits a density $f$ with respect to the Lebesgue measure, no assumption is made on the outliers $X_{i \in \mathcal{O}}$. Hence, these outlying points can be dependent, adversarial, or not even drawn from a proper probability distribution. \\

\noindent
The $\OUI$ framework is related to the well-known Huber's $\varepsilon$-contamination model \citep{huber1992robust} where it is assumed that data are i.i.d. with distribution $g = \varepsilon f_{\mathcal{I}} + (1 - \varepsilon) f_{\mathcal{O}}$, and $\varepsilon \in [0, 1)$; the distribution $f_{\mathcal{I}}$ being related to the inliers and $f_{\mathcal{O}}$ to the outliers.
However, there are several important differences. First, in the $\OUI$ the proportion of outliers is fixed and equals to $\lvert \mathcal{O} \rvert /n$, whereas it is random in the Huber's $\varepsilon$-contamination model \citep{lerasle2019lecture}.
Second, the $\OUI$ is less restrictive. Indeed, contrary to Huber's model which considers that inliers and outliers are respectively i.i.d from the same distributions, $\OUI$ does not make a single assumption on the outliers.

\subsection{MoM-KDE}
We now present our main contribution, a robust kernel density estimator based on the MoM. This estimator is essentially motivated by the fact that the classical kernel density estimation at one point corresponds to an empirical average (see Equation \eqref{eq:kde}). Therefore, the MoM principle appears to be an intuitive solution to build a robust version of the KDE. A formal definition of MoM-KDE is given below.

\begin{Definition} (MoM Kernel Density Estimator)
Let $1 \leq S \leq n$, and let $B_1, \cdots , B_S$ be a random partition of $\{1, \cdots, n\}$ into $S$ non-overlapping blocks $B_s$ of equal size $n_s \triangleq n/S$.
The MoM Kernel Density Estimator (MoM-KDE) of $f$ at $x_0$  is given by
\begin{equation}
    \hat{f}_{MoM}(x_0) \propto \emph{\text{Median}}\left(\hat{f}_{n_1}(x_0), \cdots, \hat{f}_{n_S}(x_0) \right) \; ,
\end{equation}
where $\hat{f}_{n_s}(x_0)$ is the value of the standard kernel density estimator at $x_0$ obtained via the samples of the $s$-th block $B_s$.
Note that $\hat{f}_{MoM}(\cdot)$ is not necessarily a density as its integral may not be equal to $1$. When needed, we thus normalize it by its integral the same way it is proposed by \citet{devroye2012combinatorial}.
\end{Definition}

\begin{figure}[t]
\centering
\begin{minipage}[t]{0.47\linewidth}
    \centering
    \includegraphics[width=\linewidth]{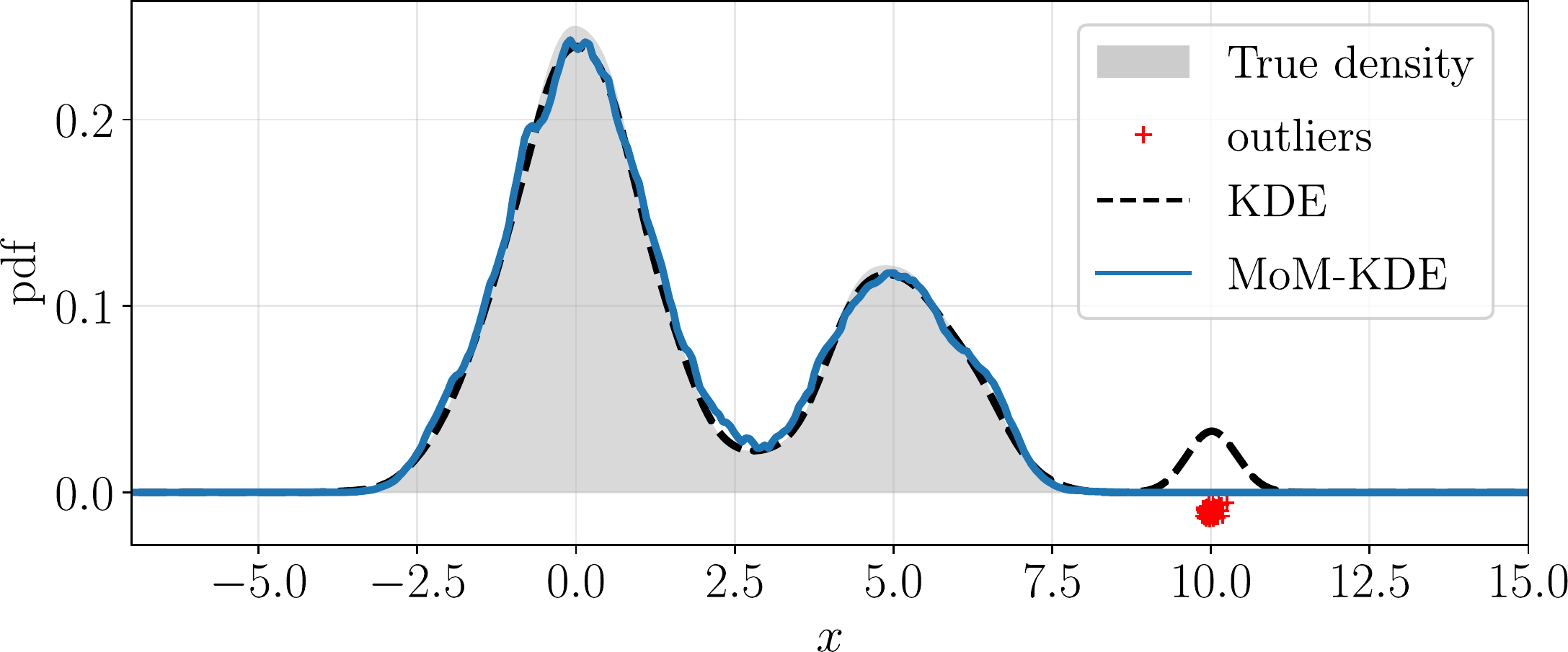}\\
    {\small (a)\; One-dimensional}
\end{minipage}\hfill
\begin{minipage}[t]{0.53\linewidth}
    \centering
    \includegraphics[width=\linewidth]{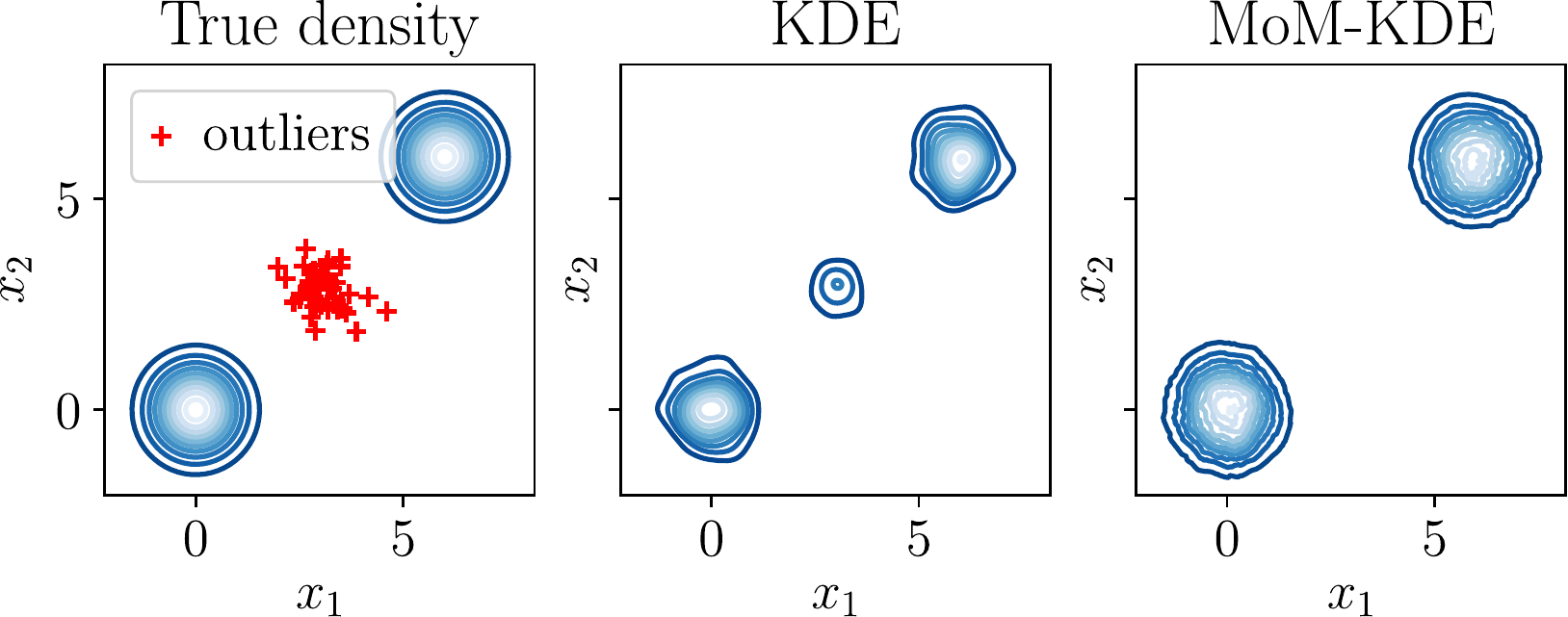}\\
    {\small (b)\; Two-dimensional}
\end{minipage}
\caption{True density, outliers, KDE, and MoM-KDE. (a) Estimates from a $1$-D true density and outliers from a normal density centered in $\mu_{\calO}=10$ with variance $\sigma^2_{\calO}=0.1$. (b) Estimates from a $2$-D true density and outliers from a normal density centered in $\mu_{\calO}=(3, 3)$ with variance $\sigma^2_{\calO}=0.5 I_2$.}
\label{fig:bump_kde}
\end{figure}

Broadly speaking, MoM estimators appear to be a good tradeoff between the unbiased but non robust empirical mean and the robust but biased median \citep{lecue2018robust}. A visual example of the robustness of MoM-KDE is displayed in Figure \ref{fig:bump_kde}.
We now give a simple example highlighting the robustness of MoM-KDE.
\begin{Example}(MoM-KDE v.s.\ Uniform KDE) Let the inliers be i.i.d. samples from a uniform distribution on the interval $[-1,1]$ and the outliers be i.i.d. samples from another uniform distribution on $[-3,3]$. Let the kernel function be the uniform kernel, $x_0=2$ and $h\in (0,1)$. Then if $S>2|\calO|$, we obtain
\begin{equation*}
    \lvert\hat{f}_{MoM}(x_0) - f(x_0)\rvert = 0 \quad  a.s. \quad \text{and} \quad \IP\left(\lvert\hat{f}_{n}(x_0) - f(x_0)\rvert = 0\right) = (1-h/3)^{|\calO|}\neq 1 \; .
\end{equation*}
\end{Example}
This result shows that the MoM-KDE makes (almost surely) no error at the point $x_0$. On the contrary, the KDE here has a non-negligible probability to make an error.

\subsection{Time complexity}
The complexity of MoM-KDE to evaluate one point is the same as the standard KDE, $\calO(n)$; $\calO(S\cdot \frac{n}{S})$ for the block-wise evaluation and $\calO(n)$ to compute the median with the \textit{median-of-medians algorithm} \citep{blum1973time}.
Since RKDE and SPKDE are KDEs with modified weights, they also perform the evaluation step in $\calO(n)$ time. However, these weights need to be learnt, thus requiring an additional non-negligible computing capacity.
Indeed, each one of them rely on an iterative method -- the iteratively reweighted least squares algorithm and the projected gradient descent algorithm, that both have a complexity of $\calO(n_{iter} \cdot n^2)$, where $n_{iter}$ is the number of needed iterations to reach a reasonable accuracy. MoM-KDE on the other hand does not require any learning procedure. Note that the evaluation step can be accelerated through several ways, hence potentially reducing computational time of all these competing methods \citep{gray2003nonparametric, gray2003rapid, wang2019nonparametric, backurs2019space}. Theoretical time complexities are gathered in Table \ref{tab:complexity}.

\begin{table}
	\footnotesize
	\centering
	\ra{1.3}
\caption{Computational complexity}
\label{tab:complexity}
\begin{adjustbox}{max width=\textwidth}
	\begin{tabular}{@{}llll@{}}
		\toprule
		\textbf{Method} & \textbf{Learning} & \textbf{Evaluation} & \textbf{Iterative method} \\ \midrule
		\textbf{KDE} \citep{parzen1962estimation} & \qquad -- &  $\calO(n)$ & no \\
		\textbf{RKDE} \citep{kim2012robust} & $\calO(n_{iter} \cdot n^2)$ & $\calO(n)$ & yes \\
		\textbf{SPKDE} \citep{vandermeulen2014robust}& $\calO(n_{iter} \cdot n^2)$ & $\calO(n)$ & yes\\
		\textbf{MoM-KDE} & \qquad -- & $\calO(n)$ & no\\
		\bottomrule
	\end{tabular}
\end{adjustbox}
\medskip
\end{table}

\section{Theoretical analysis}
In this section, we give a finite-sample high-probability error bound in the $L_\infty$-norm for MoM-KDE under the $\OUI$ framework. To our knowledge, we are the first to provide such error bounds in robust kernel density estimation under this framework. In particular, our objective is to prove that even with a contaminated dataset, MoM-KDE achieves a similar convergence rate than KDE without outliers \citep{sriperumbudur2012consistency, jiang2017uniform, wang2019dbscan}. In order to build this high-probability error bound, it is assumed, among other standard hypotheses, that the true density is Hölder-continuous, a smoothness property usually considered in KDE analysis \citep{tsybakov2008introduction, jiang2017uniform, wang2019dbscan}. In addition, we show the consistency in the $L_1$-norm. In this last result, we will see that the aforementioned assumptions are not necessary to obtain the consistency. In the following, we give the necessary definitions and assumptions to perform our non-asymptotic analysis.

\subsection{Setup and assumptions}

Let us first list the usual assumptions, notably on the considered kernel function, that will allow us to derive our results. They are standard in KDE analysis, and are chosen for their simplicity of comprehension \citep{tsybakov2008introduction, jiang2017uniform}. More general hypotheses could be made in order to obtain the same results, notably assuming kernel of order $\ell$ (see for example the works of  \citet{tsybakov2008introduction} and \citet{wang2019dbscan}).

\begin{itemize}[leftmargin=*]
    \item[] \textbf{Assumption 1.} (Bounded density) $\|f\|_{\infty}<\infty$. \vspace{.5em}
    \item[] We make the following assumptions on the kernel $K$.
    \item[] \textbf{Assumption 2.} (Density kernel) $\forall u \in \IR^d, K(u)\geq 0$, and $\displaystyle\int K(u)du=1$.
    \item[] \textbf{Assumption 3.} (Spherically symmetric and non-increasing) There exists a non-increasing function $k : \IR_{+} \longrightarrow \IR_{+}$ such that $K(u) = k(\|u\|)$ for all $u \in \IR^d$, where $\|\cdot\|$ is any norm of $\IR^d$.
    \item[] \textbf{Assumption 4.} (Exponentially decaying tail)
There exists positive constants $\rho$, $C_{\rho}, t_0 > 0$ such that for all $ t>t_0$ $$k(t) \leq C_{\rho} \cdot \exp(-t^\rho)  \; .$$
\end{itemize}

All the above assumptions are respected by most of the popular kernels, in particular the Gaussian, Exponential, Uniform, Triangular, Cosine kernel, etc. 
Furthermore, the last assumption implies that for any $m>0$, we have $\int \|u\|^mK(u)du < \infty$ (finite norm moment) \citep{jiang2017uniform}.
Finally, when taken together, these assumptions imply that the kernel satisfies the VC property \citep{wang2019dbscan}.
Theses are key properties to provide the bounds presented in the next section. \\

\noindent
Before stating our main results, we recall the definition of the Hölder class of functions.
\begin{Definition} \label{def:holder} (Hölder class) Let $T$ be an interval of $\IR^d$, and $0<\alpha \leq 1$ and $L>0$ be two constants. We say that a function $f:T\rightarrow \IR$ belongs to the Hölder class $\Sigma(L, \alpha)$ if it satisfies
\begin{equation}
    \forall x,x' \in T, \quad |f(x)-f(x')| \leq L\|x-x'\|^{\alpha} \; .
\end{equation}
\end{Definition}

This definition implies a smoothness regularization on the function $f$, and is a convenient property to bound the bias of KDE-based estimators.

\subsection{{\boldmath $L_\infty$} and {\boldmath${L_1}$} consistencies of MoM-KDE}
This section states our central finding, a $L_\infty$ finite-sample error bound for MoM-KDE that proves its consistency and yields the same convergence rate as KDE with uncontaminated data. The latter is given by the following Lemma partly proven by \citet{sriperumbudur2012consistency} and verified several times in the literature \citep{gine2002rates, jiang2017uniform, wang2019dbscan}.

\begin{lemma}\label{lemma:kde}($L_\infty$ error-bound of the KDE without anomalies) Suppose that $f$ belongs to the class of densities $\mathcal{P}(\alpha, L)$ 
defined as
\begin{equation}
    \mathcal{P}(\alpha, L) \triangleq \left\{ f \mid f \geq 0, \int f(x)dx = 1, \text{ and } f \in \Sigma(\alpha, L) \right\} \;,
\end{equation}
where $\Sigma(\alpha, L)$ is the Hölder class of function on $\IR^d$ (Definition \ref{def:holder}). 
Grant assumptions $1$ to $4$ and let $n>1$, $h\in (0,1)$ and $S\geq1$ such that $nh^d\geq S$ and $nh^d\geq\lvert \log(h) \rvert$. Then with probability at least $1-\exp(-S)$, we have 
\begin{equation}
    \label{eq:bound_kde}
    \|\hat{f}_{n} - f\|_\infty \leq C_1\sqrt{\frac{S|\log(h)|}{nh^d}} + C_2h^\alpha \;,
\end{equation}
where $C_2=L\displaystyle\int\|u\|^\alpha K(u)du<\infty$ and $C_1$ is a constant that only depends on $\|f\|_\infty$, the dimension $d$, and the kernel properties.
\end{lemma}
This Lemma comes from the well-known bias-variance decomposition, where we separately bound the variance (see Theorem $3.1$ of \citet{sriperumbudur2012consistency}) and the bias (see e.g. \citep{tsybakov2008introduction} or \citep{rigollet2009optimal}). It shows the consistency of KDE without anomalies, as soon as $h\rightarrow 0$ and $nh^d\rightarrow \infty$. \\[0.1cm]

\noindent
We now present our main result. Its objective is to show that even under the $\OUI$ framework, we do not need any additional hypothesis -- besides the ones of the previous lemma -- to show that MoM-KDE achieves the same convergence rate as KDE when used with uncontaminated data. \\

\begin{proposition} ($L_\infty$ error-bound of the MoM-KDE under the $\OUI$)
Suppose that $f$ belongs to the class of densities $\mathcal{P}(\alpha, L)$ and grant assumptions $1$ to $4$. Let $S$ be the number of blocks, $\delta>0$ such that $S>(2+\delta)|\mathcal{O}|$, and $\Delta = (1/(2+\delta) - |\mathcal{O}|/S)$.
Then, for any  $h\in(0,1)$, $\delta$ sufficiently small, and $n \geq 1$ such that $nh^d \geq S\log(2(2+\delta)/\delta)$, and $nh^d  \geq S|\log(h)|$, we have with probability at least $1-\exp(-2\Delta^2S)$,
\begin{equation}
    \label{eq:bound_mom}
    \|\hat{f}_{MoM} - f\|_\infty \leq C_1\sqrt{\frac{S \log\left(\frac{2(2+\delta)}{\delta}\right)|\log(h)|}{nh^d}} + C_2h^\alpha \;,
\end{equation}
where $C_2=L\displaystyle\int\|u\|^\alpha K(u)du<\infty$ and $C_1$ is a constant that only depends on $\|f\|_\infty$, the dimension $d$, and the kernel properties.
\end{proposition}

The proof is given in the supplementary material. From equation (\ref{eq:bound_mom}), the optimal choice of the bandwidth is $h \asymp \left(\dfrac{S \log(n)}{n}\right)^{1/(2\alpha + d)}$ leading to the final rate of $\left(\dfrac{S \log(n)}{n}\right)^{\alpha/(2\alpha + d)}.$ This convergence rate is the same (up to a constant) to the one of KDE without anomalies, with the same exponential control (Lemma \ref{lemma:kde}). Note that when there is no outlier, i.e. $|\mathcal{O}| = 0$, the bound holds for $S=1$, and we recover the classical KDE minimax optimal rate \citep{wang2019dbscan}. In addition, the previous proposition states that the convergence of the MoM-KDE only depends on the number of outliers in the dataset, and not on their ``type''. This estimator is therefore robust in a wide range of scenarios, including the adversarial one. \\

\noindent
We now give a $L_1$-consistency result under mild hypotheses, which is known to reflect the global performance of the estimate. Indeed, small $L_1$ error leads to accurate probability estimation \citep{devroyegy}.
\begin{proposition}($L_1$-consistency in probability) If $n/S\rightarrow \infty$, $h\rightarrow 0$, $n h^d \rightarrow \infty$, and $S>2|\mathcal{O}|$, then
\begin{equation}
  \|\hat{f}_{MoM} - f\|_1 \overset{\mathcal{P}}{\underset{n\rightarrow \infty}{\longrightarrow}} 0 \; .
\end{equation}
\end{proposition}
This result is obtained by bounding the left-hand part by the errors in the healthy blocks only, i.e. those without anomalies. Under the hypothesis of the proposition, these errors are known to converge to $0$ in probability \citep{wang2019dbscan}. The complete proof is given in supplementary material.
Contrary to SPKDE \citep{vandermeulen2014robust}, no assumption on the outliers generation process is necessary to obtain this consistency result. Moreover, while we need to assume that the proportion of outliers is perfectly known to prove the convergence of SPKDE, the MoM-KDE converges whenever the number of outliers is overestimated.

\section{Numerical experiments}
In this section, we display numerical results supporting the relevance of MoM-KDE. All experiences were run over a personal laptop computer using \texttt{Python}. The code of MoM-KDE is made available online\footnote{\url{https://github.com/lminvielle/mom-kde} \quad For the sake of comparison, we also implemented RKDE and SPKDE.}.

\paragraph{Comparative methods.}
In the following experiments, we propose to compare MoM-KDE to the classical KDE and two robust versions of KDE, called RKDE \citep{kim2012robust} and SPKDE \citep{vandermeulen2014robust}.

As previously explained, RKDE takes the ideas of robust M-estimation and translate it to kernel density estimation. Authors point out that classical KDE estimator can be seen as the minimizer of a squared error loss in the Reproducing Kernel Hilbert Space $\mathcal{H}$ corresponding to the chosen kernel. Instead of minimizing this loss, they propose to minimize a robust version of it, $\sum_i\rho(\|\phi(X_i)-g\|_\mathcal{H})$, with respect to $g\in \mathcal{H}$. Here $\phi$ is the canonical feature map and $\rho(\cdot)$ is either the robust Huber or Hampel function. The solution of the newly expressed problem is then found using the iteratively reweighted least squares algorithm.

SPKDE proposes to scale the standard KDE in a way that it decontaminates the dataset. This is done by minimizing the function $\|\beta \hat{f}_n - g\|_2$ with respect to $g$, belonging to the convex hull of $\{k_h(\cdot, X_i)\}_{i=1}^n$. Here, $\beta$ is an hyperparameter that controls the robustness and $\hat{f}_n$ is the KDE estimator. The minimization is shown to be equivalent to a quadratic program over the simplex, solved via projected gradient descent.

\paragraph{Metrics.}
The performance of the MoM-KDE is measured through three metrics, two are used to measure the similarity between the estimated and the true density, and one describes performances of an anomaly detector based on the estimated density.
The first one is the Kullback-Leibler divergence \citep{kullback1951information} which is the most used in robust KDE \citep{kim2008robust, kim2011robustness, kim2012robust, vandermeulen2014robust}. Used to measure the similarity between distributions, it is defined as $$D_{\text{KL}}(\hat{f} \lVert f) = \displaystyle\int \hat{f}(x) \log\left(\dfrac{\hat{f}(x)}{f(x)}\right) dx \; .$$

As the Kullback-Leibler divergence is non-symmetric and may have infinite values when distributions do not share the same support, we also consider the Jensen-Shannon divergence \citep{endres2003new, liese2006divergences}.
It is a symmetrized version of $D_{\text{KL}}$, with positive values, bounded by $1$ (when the logarithm is used in base $2$), and has found applications in many fields, such as deep learning~\citep{goodfellow2014generative} or transfer learning~\citep{segev2017learn}.
It is defined as $$D_{\text{JS}}(\hat{f} \lVert f) = \frac{1}{2} \left(D_{\text{KL}}(\hat{f} \lVert g) + D_{\text{KL}}(f \lVert g)\right), \quad \text{with} \quad g = \frac{1}{2} (\hat{f} + f )\; .$$

Motivated by real-world application, the third metric is not related to the true density, which is usually not available in practical cases.
Instead, we quantify the capacity of the learnt density to detect anomalies using the well-known Area Under the ROC Curve criterion (AUC). 
An input point $x_0$ is considered abnormal whenever $\hat{f}(x_0)$ is below a given threshold.

\paragraph{Hyperparameters.}
All estimators are built using the Gaussian kernel. The number of blocks in MoM-KDE is selected on a regular grid of $20$ values between $1$ and $2 \vert \calO \rvert + 1$ in order to obtain the lowest $D_{\text{JS}}$. The bandwidth $h$ is chosen for KDE via the pseudo-likelihood $k$-cross-validation method~\citep{turlach1993bandwidth}, and used for all estimators. The construction of RKDE follows exactly the indications of its authors \citep{kim2012robust} and $\rho(\cdot)$ is taken to be the Hampel function as they empirically showed that it is the most robust. For SPKDE, the true ratio of anomalies is given as an input parameter.

\subsection{Results on synthetic data.}
\newcommand{\ratio}{0.24}
\newcommand{\ratiob}{0.8}
\begin{figure}
    \centering
    \begin{minipage}[t]{\ratio\linewidth}
        \centering
        \includegraphics[height=\ratiob\linewidth]{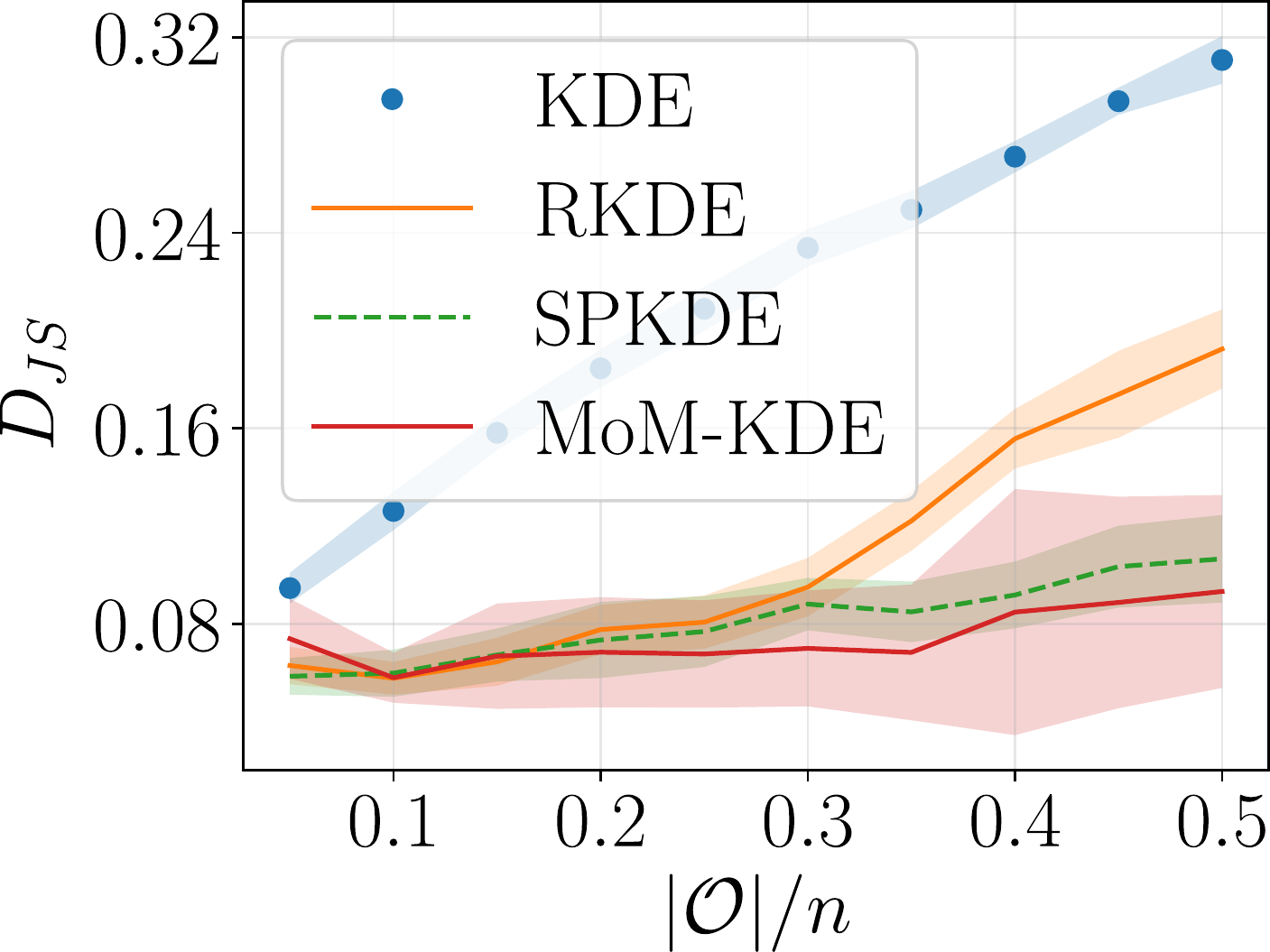}
        {\small (a)\; Uniform}
    \end{minipage}\hfill
    \begin{minipage}[t]{\ratio\linewidth}
        \centering
        \includegraphics[height=\ratiob\linewidth]{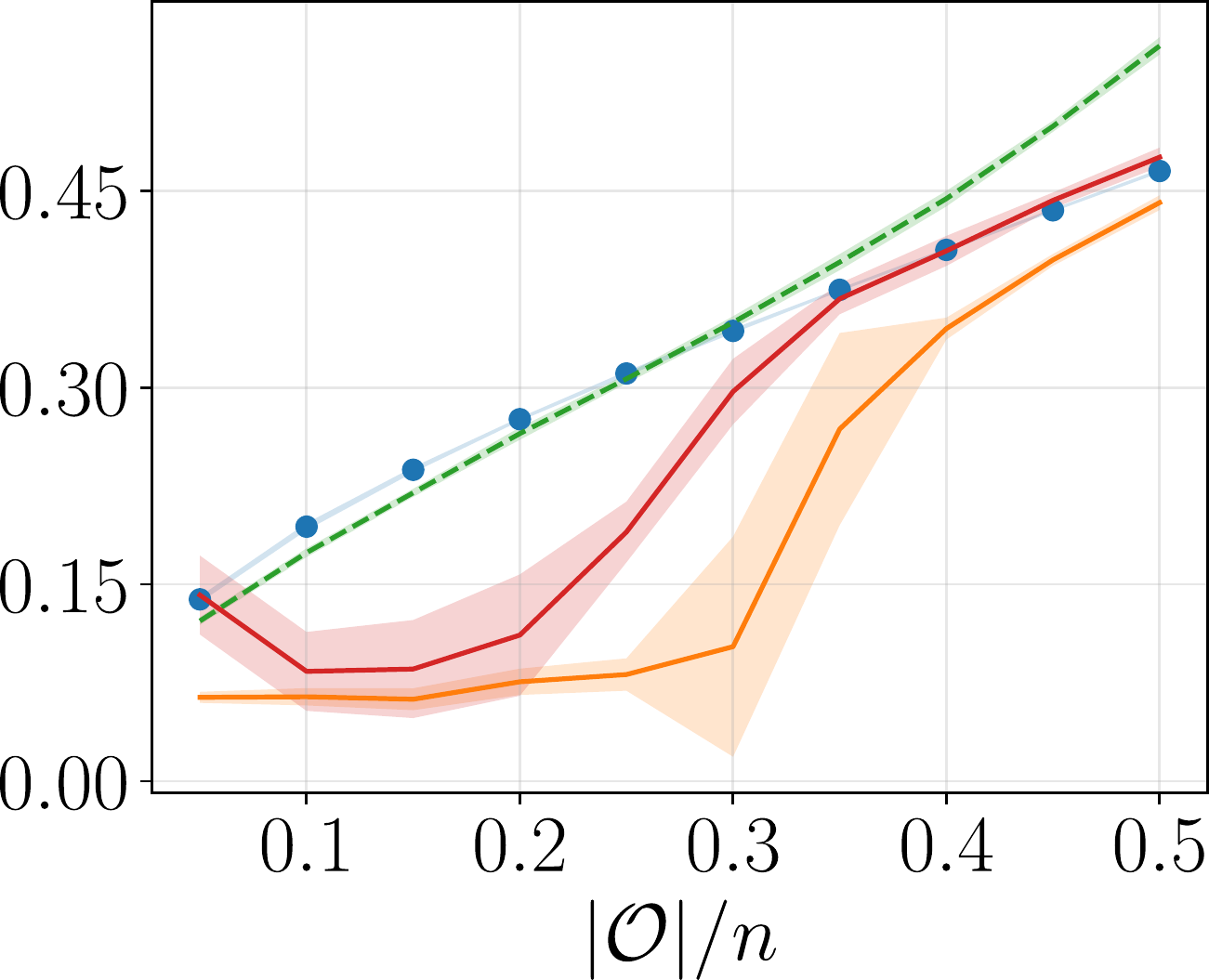}
        {\small (b)\; Regular Gaussian}
    \end{minipage}
    \begin{minipage}[t]{\ratio\linewidth}
        \centering
        \includegraphics[height=0.81\linewidth]{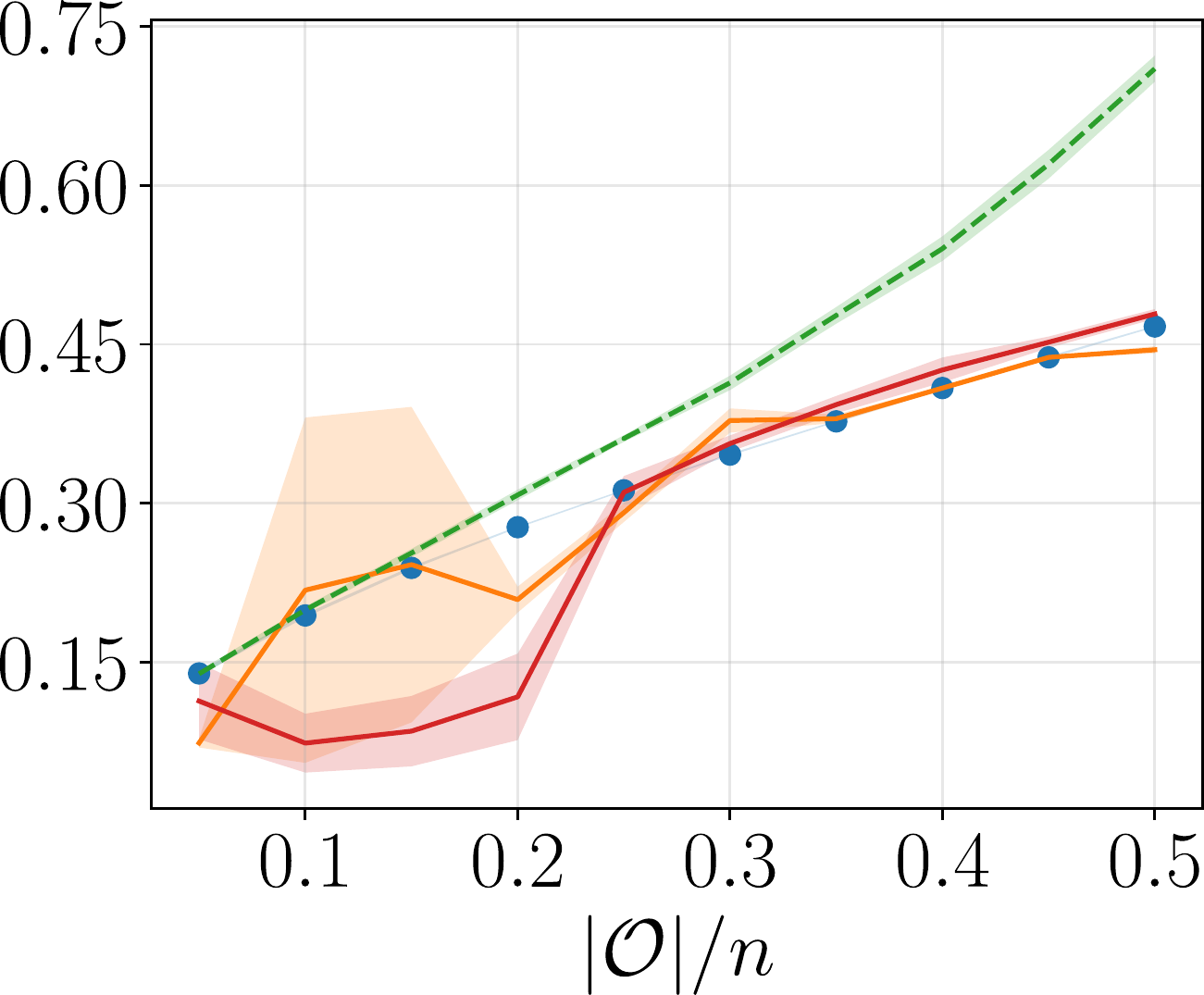}
        {\small (c)\; Thin Gaussian}
    \end{minipage}
    \begin{minipage}[t]{\ratio\linewidth}
    \centering
    \includegraphics[height=\ratiob\linewidth]{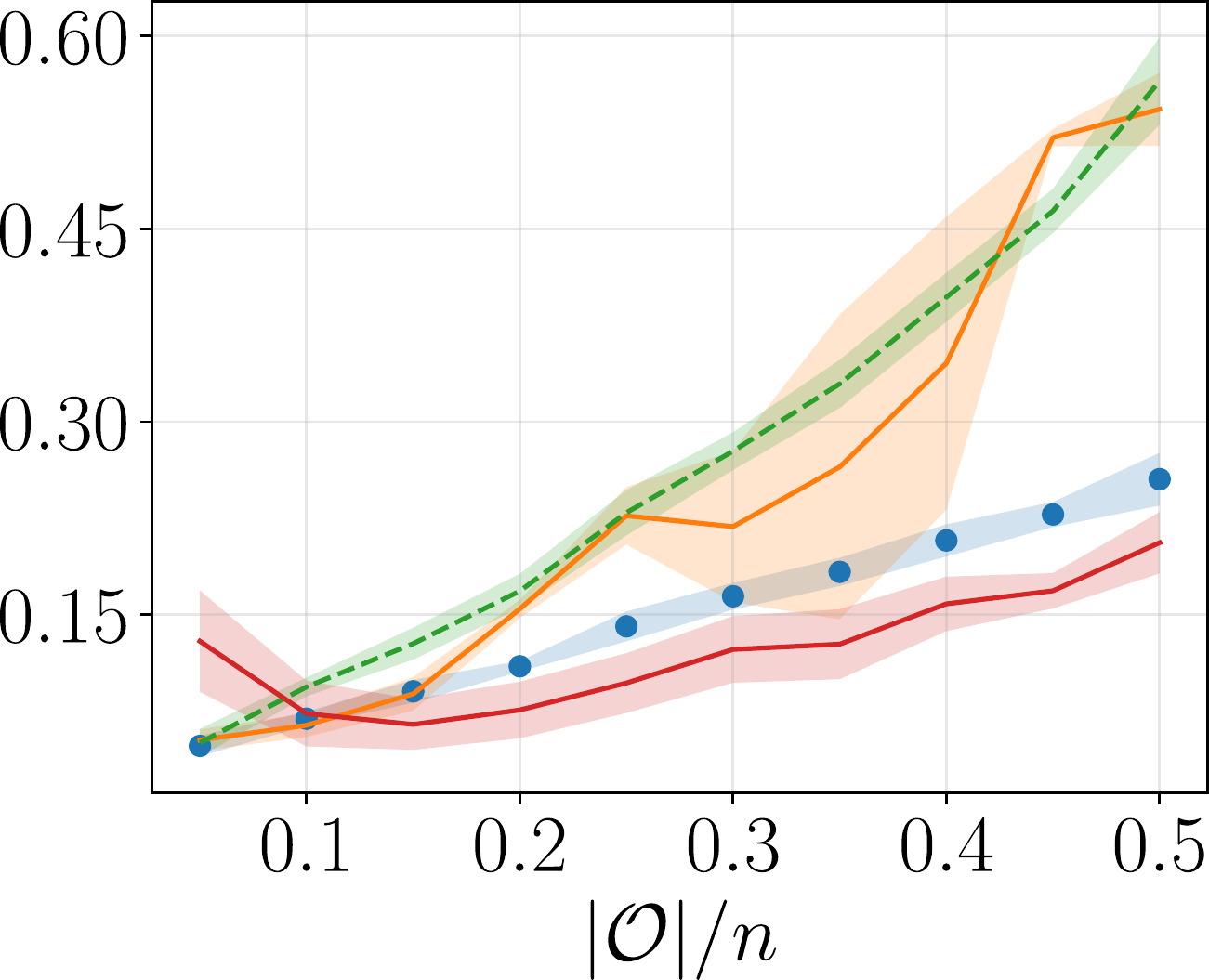}
    {\small (d)\; Advsersarial Thin Gaussian}
    \end{minipage}
    \caption{Density estimation with synthetic data. The displayed metric is the Jensen-Shannon divergence. A lower score means a better estimation of the true density.}
    \label{fig:experiments_synth}
\end{figure}
To evaluate the efficiency of the MoM-KDE against KDE and its robust competitors, we set up several outlier situations. In all theses situations, we draw $N=1000$ inliers from an equally distributed mixture of two normal distribution $\mathcal{N}(\mu_1, \sigma_1 )$ and $\mathcal{N}(\mu_2, \sigma_2)$ with $\mu_1 = 0$, $\mu_2=6$, and $\sigma_1 = \sigma_2 = 0.5$.
The outliers however are sampled through various schemes:

\begin{enumerate}[leftmargin=1cm]
    \item [(a)]\textbf{Uniform.} A uniform distribution $U([\mu_1-3, \mu_2+3])$ which is the classical setting used for outlier simulation.
    \item [(b)] \textbf{Regular Gaussian.} A \emph{similar}-variance normal distribution $\mathcal{N}(3, 0.5)$ located between the two inlier clusters.
    \item [(c)] \textbf{Thin Gaussian.} A \emph{low}-variance normal distribution $\mathcal{N}(3,0.01)$ located between the two inliers clusters.
    \item [(d)] \textbf{Adversarial Thin Gaussian.} A low variance normal distribution $\mathcal{N}(0, 0.01)$ located on one of the inliers' Gaussian mode.
    This scenario can be seen as adversarial as an ill-intentioned agent may hide wrong points in region of high density. It is the most challenging setting for standard robust estimators as they are in general robust to outliers located outside the support of the density we wish to estimate.
\end{enumerate}

For all situations, we  consider several ratios of contamination and set the number of outliers $\lvert \calO \rvert$ in order to obtain a ratio $\lvert \calO \rvert/n$ ranging from $0.05$ to $0.5$ with $0.05$-wide steps. Finally, to evaluate the pertinence of our results, for each set of parameters, data are generated $10$ times.\\

\noindent
We display in Figure \ref{fig:experiments_synth} the results over synthetic data using the $D_{\text{JS}}$ score. The average scores and standard deviations over the $10$ experiments are represented for each outlier scheme and ratio.  Overall, the results show the good performance of MoM-KDE in all the considered situations. Furthermore, they highlight the dependency of the two competitors to the type of outliers. Indeed, as SPKDE is designed to handle uniformly distributed outliers, the algorithm struggles when confronted with differently distributed outliers (see Figure \ref{fig:experiments_synth} (b, c, d)). RKDE performs generally better, but fails against adversarial contamination, which may be explained by its tendency to down-weight points located in low-density regions, which for this particular case correspond to the inliers. Results over $D_{\text{KL}}$ and AUC are reported in the supplementary materials. Generally, they show similar results and the same conclusions on the good performance of MoM-KDE can be made.

\subsection{Results on real data.}
Experiments are also conducted over six classification datasets: Banana, German, Titanic, Breast-cancer, Iris, and Digits. They contain respectively $n = 5300, 1000, 2201, 569, 150$ and $1797$ data points having
$d = 2, 20, 3, 30, 4$ and $64$ input dimensions.
They are all publicly available either from open repositories \footnote{\url{http://www.raetschlab.org/Members/raetsch/benchmark/}}
(for the first three) or directly from Scikit-learn package (for the last three)~\citep{scikit-learn}.  We follow the approach of \citet{kim2012robust} that consists in setting the class labeled $0$ as outliers and the rest as inliers.
To artificially control the outlier proportion, we randomly downsample the abnormal class to reach a ratio $\lvert \calO \rvert/n$ ranging from $0.05$ to $0.5$ with $0.05$-wide steps. When a dataset does not contain enough outliers to reach a given ratio, we similarly downsample the inliers. For each dataset and each ratio, the experiments are performed $50$ times, the random downsampling resulting in different learning datasets. The empirical performance is evaluated through the capacity of each estimator to detect anomalies, which we measure with the AUC.

Results are displayed in Figure \ref{fig:experiments_real}. With the Digits dataset, we also explore additional scenarios with changing inlier and outlier classes (specified in figure titles). Overall, results are in line with performances observed over synthetic experiments, achieving good results in comparison to its competitors. Note that even in the highest dimensional scenarios, i.e.\ Digits and Breast cancer ($d=64$ and $d=30$), MoM-KDE still behaves well, outperforming its competitors. Additional results are reported in the supplementary materials.

\renewcommand{\ratio}{0.24}
\renewcommand{\ratiob}{0.8}

\begin{figure}
    \centering
    \begin{minipage}[t]{\linewidth}
        \centering
        \begin{minipage}[t]{\ratio\linewidth}
            \centering
            \includegraphics[height=\ratiob\linewidth]{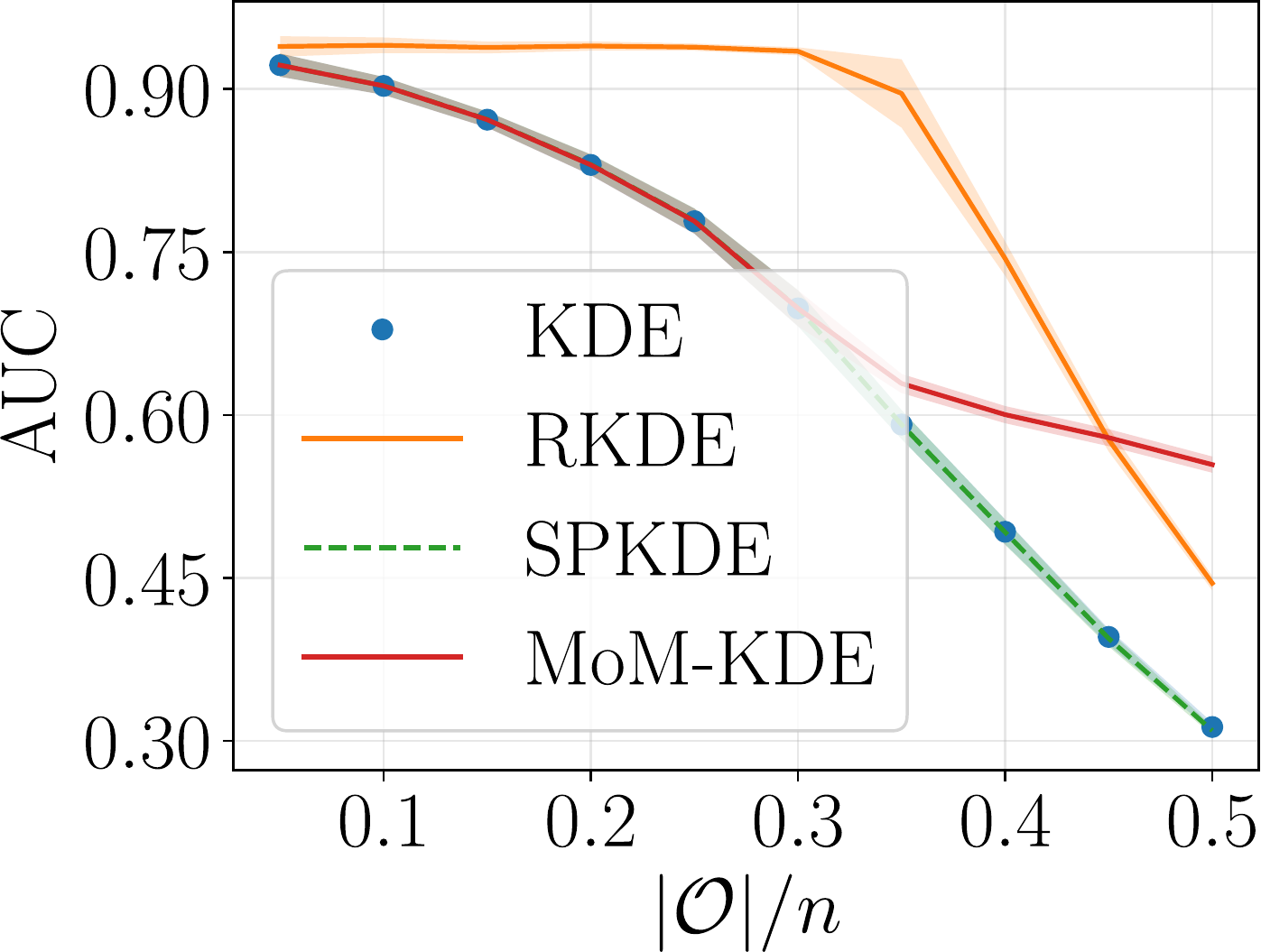}\\
            {\small (a)\; Banana}
        \end{minipage}\hfill
        \begin{minipage}[t]{\ratio\linewidth}
            \centering
            \includegraphics[height=\ratiob\linewidth]{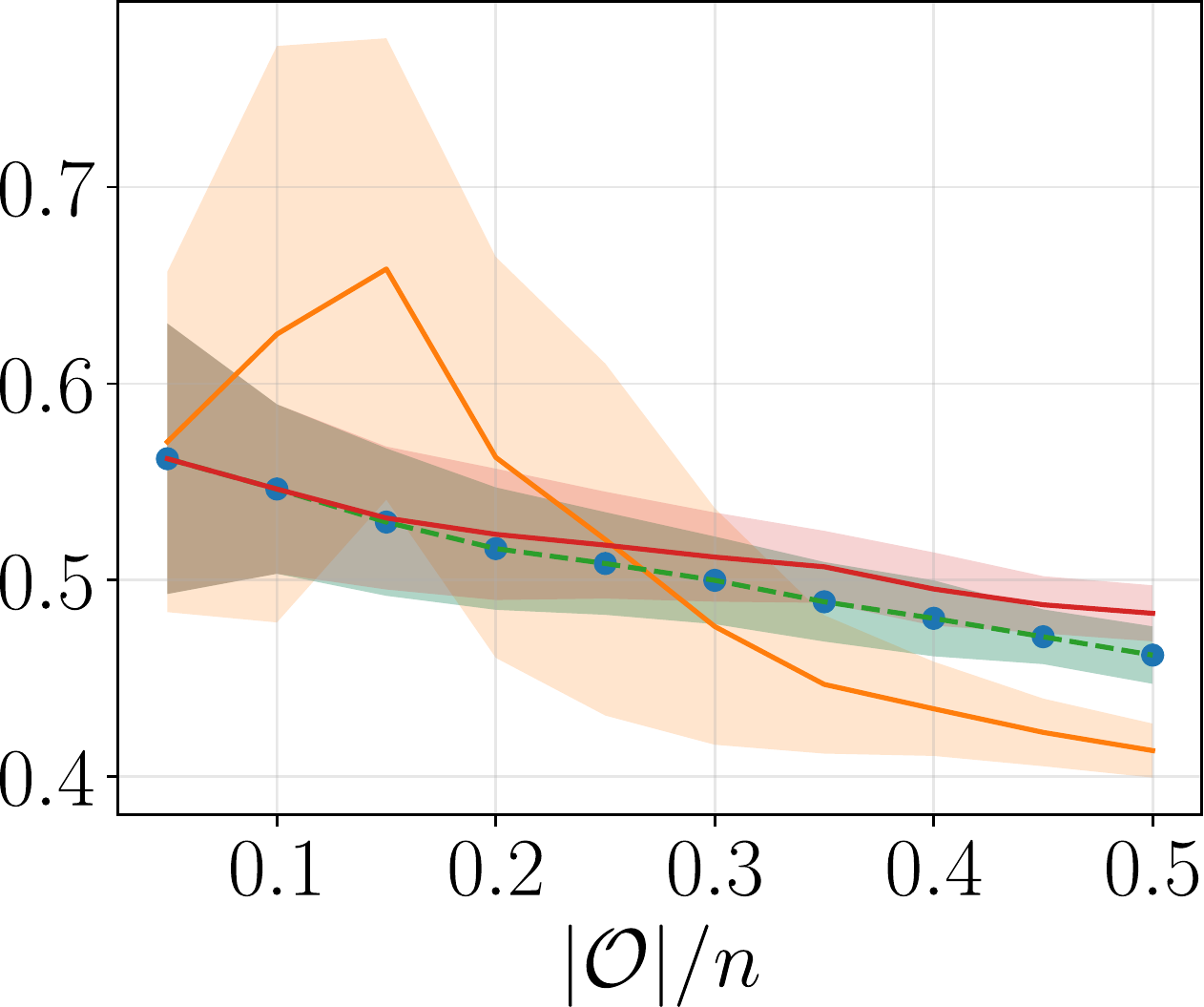}\\
            {\small (b)\; German}
        \end{minipage}
        \begin{minipage}[t]{\ratio\linewidth}
            \centering
            \includegraphics[height=\ratiob\linewidth]{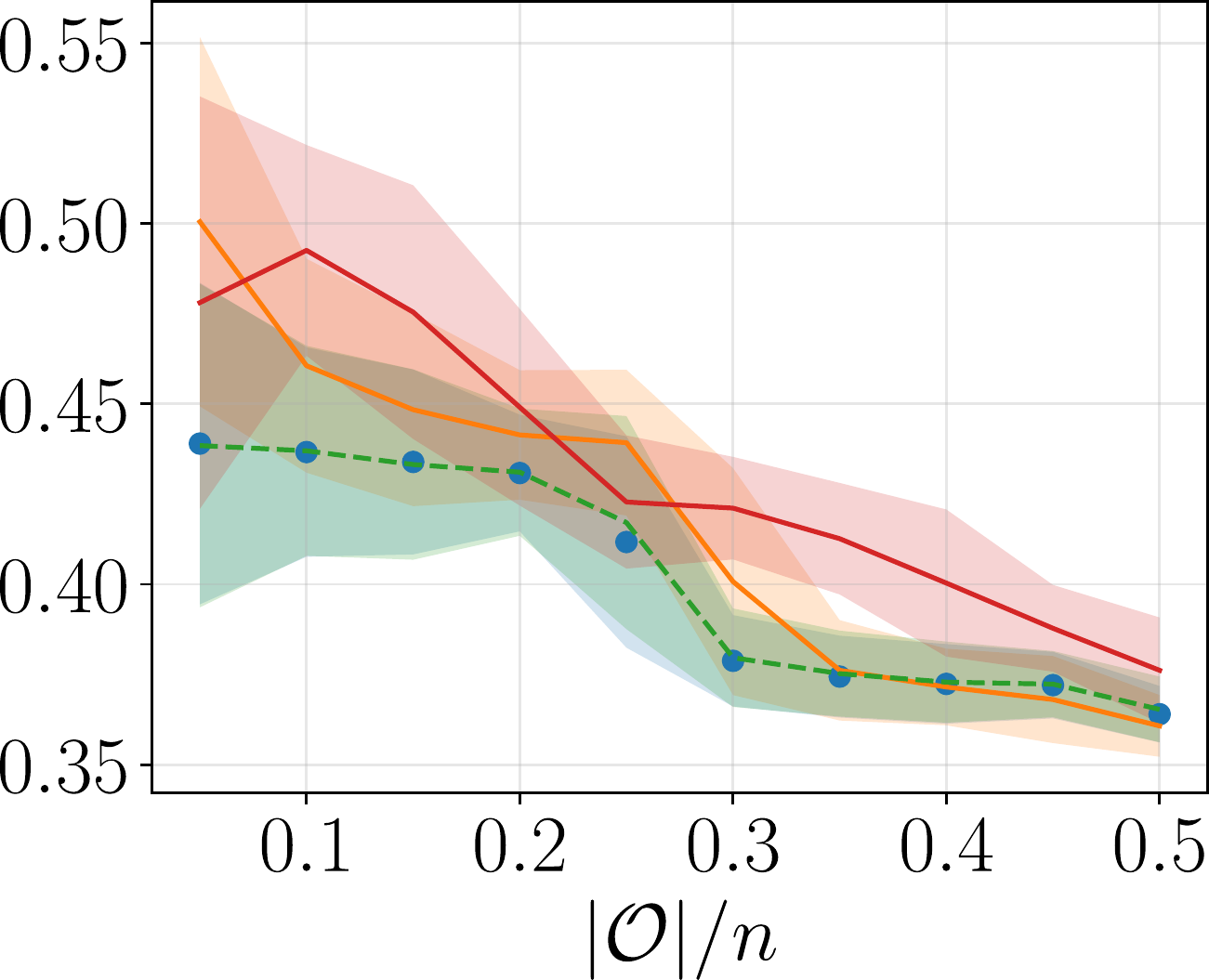}
            {\small (c)\; Titanic}
        \end{minipage}
        \begin{minipage}[t]{\ratio\linewidth}
            \centering
            \includegraphics[height=\ratiob\linewidth]{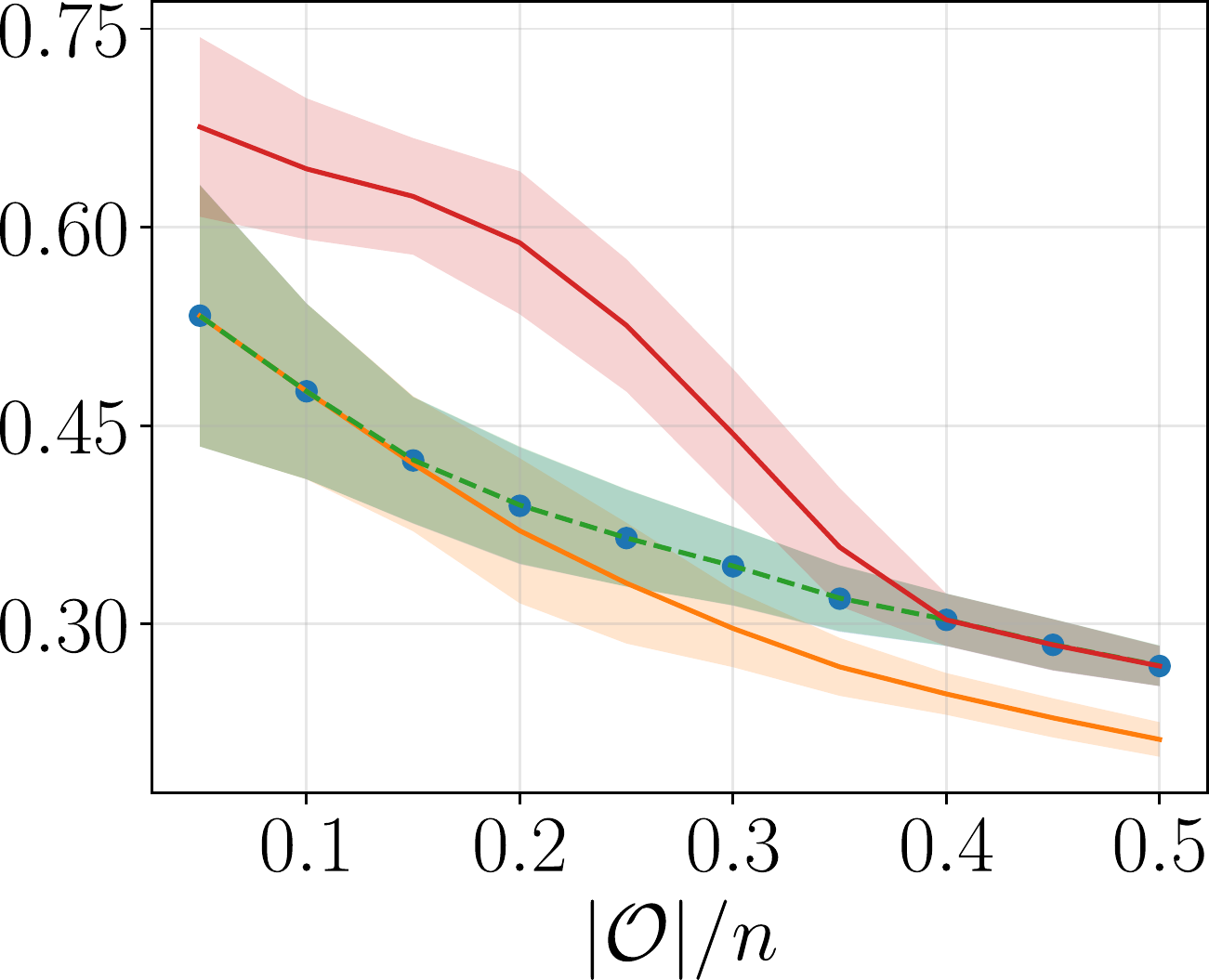}
            {\small (d)\; Breast cancer}
        \end{minipage}\\
        \bigskip
        \begin{minipage}[t]{\ratio\linewidth}
        \centering
            \includegraphics[height=\ratiob\linewidth]{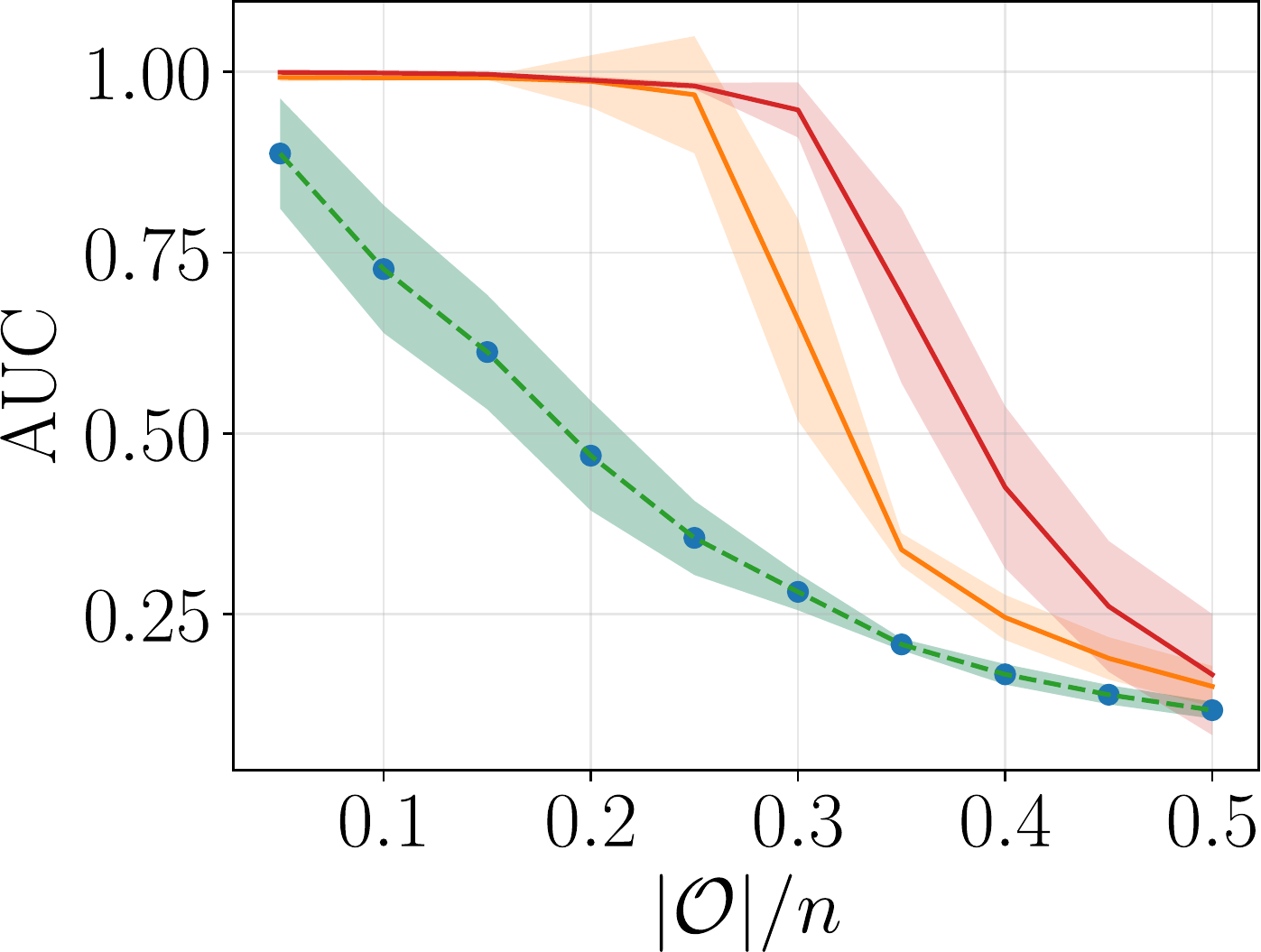}
            {\small (e)\; Iris}
        \end{minipage}\hfill
        \begin{minipage}[t]{\ratio\linewidth}
            \centering
            \includegraphics[height=\ratiob\linewidth]{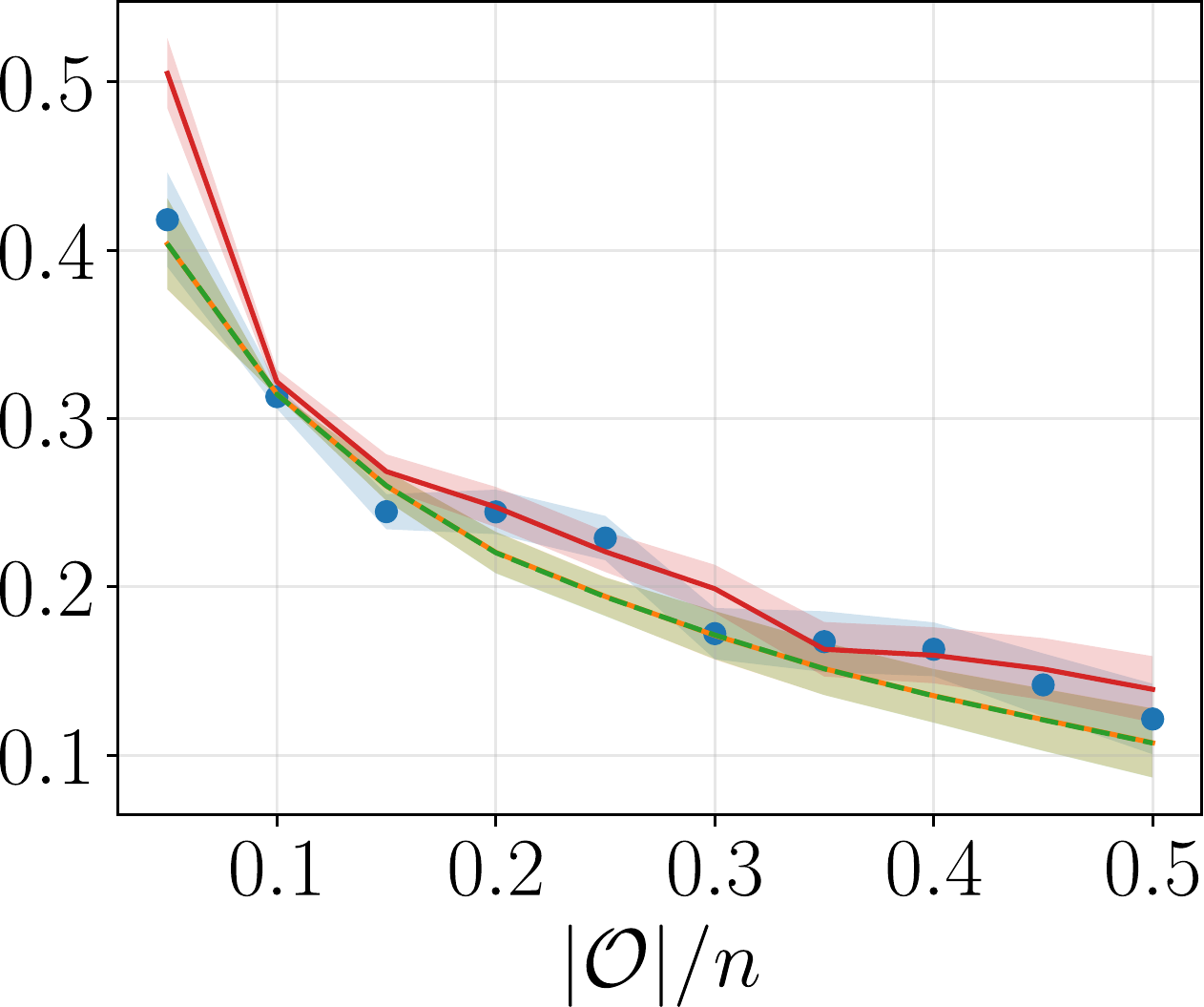}
            {\small (f)\; Digits ($\calO$: 0, $\calI$: all)}
        \end{minipage}
        \begin{minipage}[t]{\ratio\linewidth}
            \centering
            \includegraphics[height=\ratiob\linewidth]{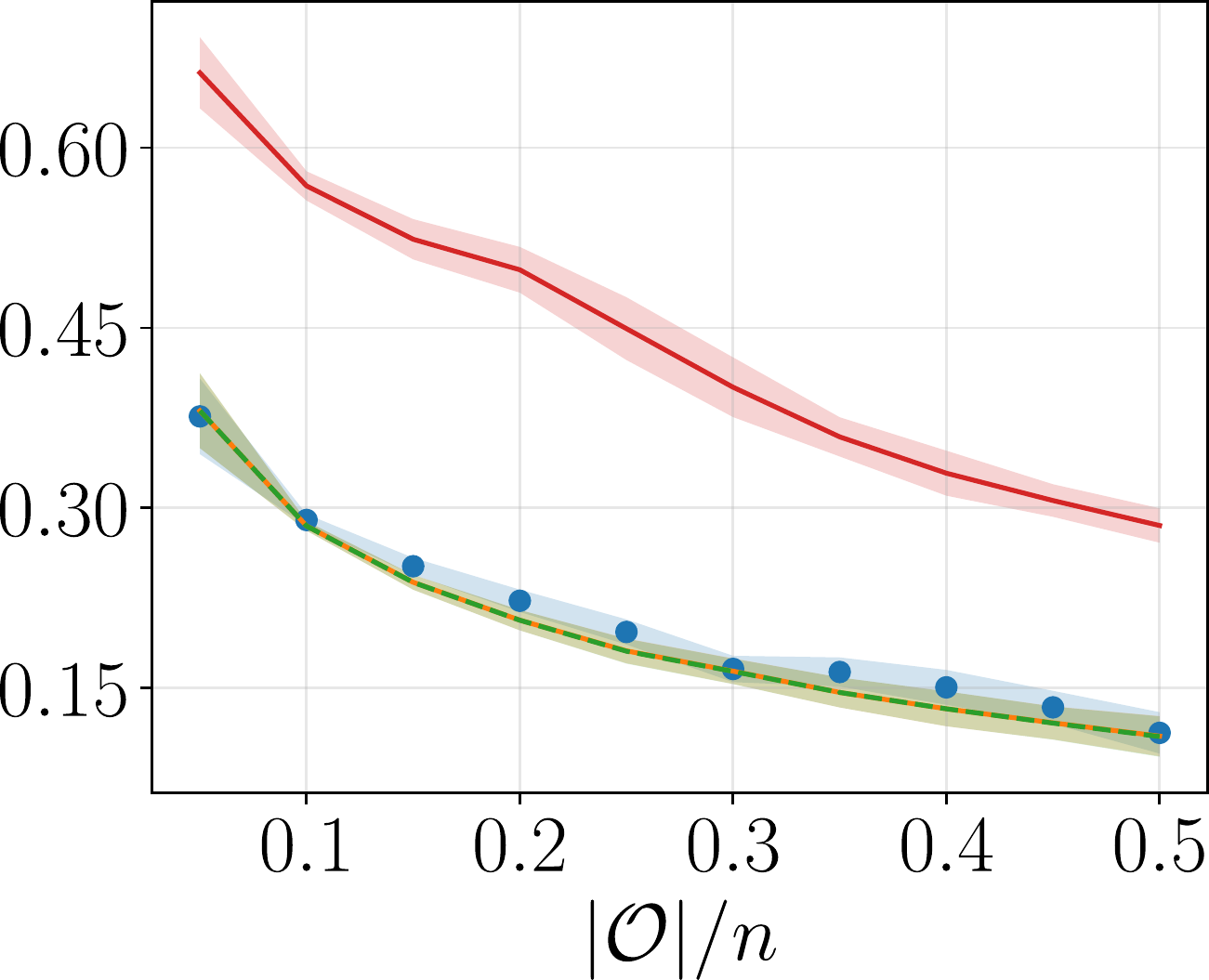}
            {\small (g)\; Digits ($\calO$: 1, $\calI$: all)}
        \end{minipage}
        \begin{minipage}[t]{\ratio\linewidth}
            \centering
            \includegraphics[height=\ratiob\linewidth]{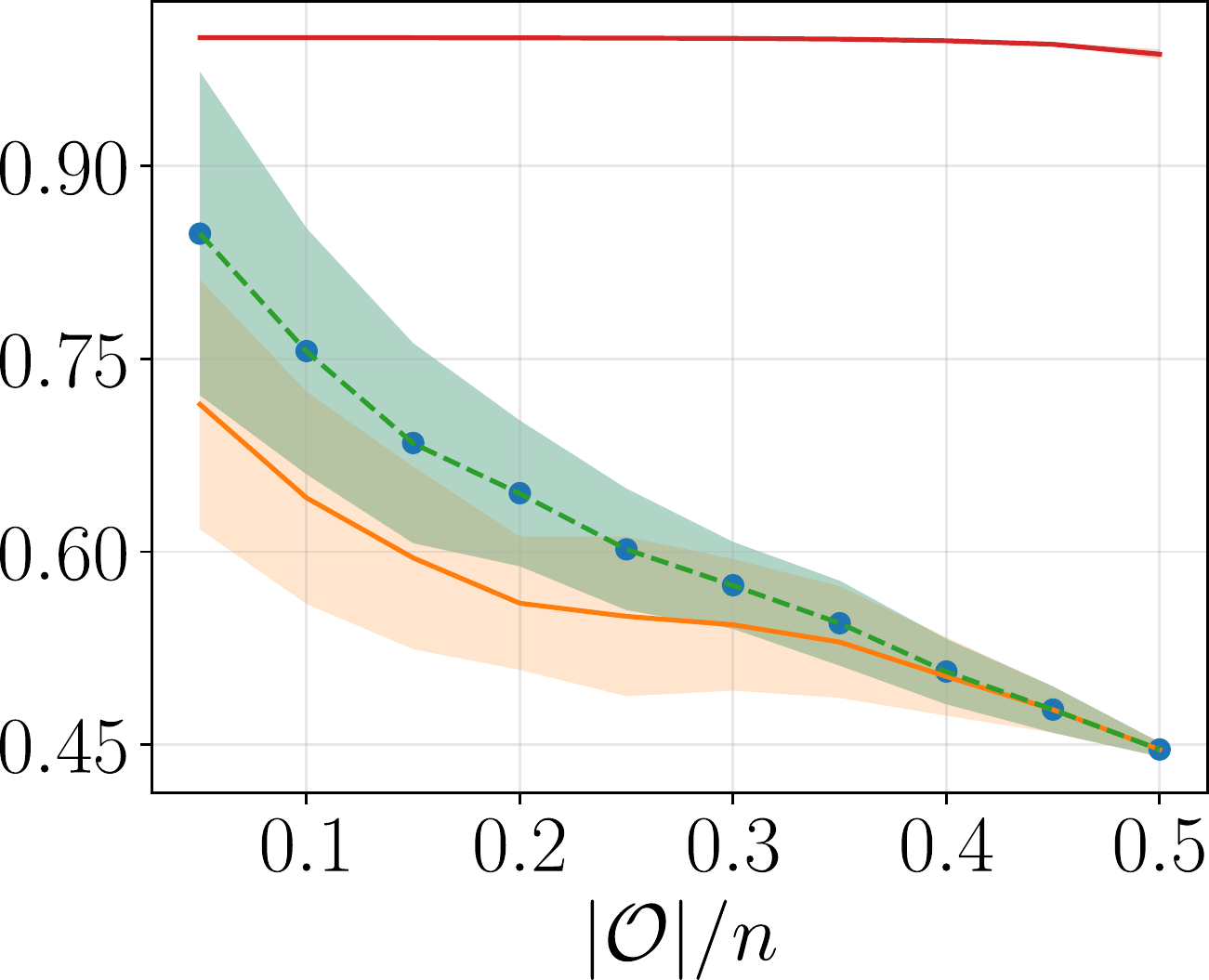}
            {\small (h)\; Digits ($\calO$: 1, $\calI$: 0)}
        \end{minipage}
    \end{minipage}
    \caption{Anomaly detection with real datasets, measured with AUC over varying outlier proportion. A higher score means a better detection of the outliers. For Digits, we specify which classes are chosen to be inliers ($\calI$) and outliers ($\calO$).}
    \label{fig:experiments_real}
\end{figure}

\section{Conclusion}
The present paper introduced MoM-KDE, a new efficient way to perform robust kernel density estimation. The method has been shown to be consistent in both $L_\infty$ and $L_1$ error-norm in presence of very generic outliers, enjoying a similar rate of convergence than the KDE without outliers. MoM-KDE achieved good empirical results in various situations while having a lower computational complexity than its competitors.

This work proposed to use the coordinate-wise median to construct its robust estimator. Future works will investigate the use of other generalization of the median in high dimension, e.g. the geometric median. In addition, further investigation will include a deeper statistical analysis under the hurdle contamination model in order to analyse the minimax optimality \citep{liu2019density} of MoM-KDE.

\balance
\bibliographystyle{abbrvnat}
\bibliography{biblio}

\appendix
\begin{center}
	\large
	APPENDIX
\end{center}

\section{Technical proofs}

\paragraph{Lemma 1.} \emph{($L_\infty$ error-bound of the KDE without anomalies) Suppose that $f$ belongs to the class of densities $\mathcal{P}(\alpha, L)$ 
defined as
\begin{equation*}
    \mathcal{P}(\alpha, L) \triangleq \left\{ f \mid f \geq 0, \int f(x)dx = 1, \text{ and } f \in \Sigma(\alpha, L) \right\} \;,
\end{equation*}
where $\Sigma(\alpha, L)$ is the Hölder class of function on $\IR^d$. 
Grant assumptions $1$ to $4$ and let $n>1$, $h\in (0,1)$ and $S\geq1$ such that $nh^d\geq S$ and $nh^d\geq\lvert \log(h) \rvert$. Then with probability at least $1-\exp(-S)$, we have 
\begin{equation*}
    \|\hat{f}_{n} - f\|_\infty \leq C_1\sqrt{\frac{S|\log(h)|}{nh^d}} + C_2h^\alpha \;,
\end{equation*}
where $C_2=L\displaystyle\int\|u\|^\alpha K(u)du<\infty$ and $C_1$ is a constant that only depends on $\|f\|_\infty$, the dimension $d$, and the kernel properties.}

\vspace{10pt}

\paragraph{Proposition 1.} \emph{($L_\infty$ error-bound of the MoM-KDE under the $\OUI$)
Suppose that $f$ belongs to the class of densities $\mathcal{P}(\alpha, L)$ and grant assumptions $1$ to $4$. Let $S$ be the number of blocks, $\delta>0$ such that $S>(2+\delta)|\mathcal{O}|$, and $\Delta = (1/(2+\delta) - |\mathcal{O}|/S)$.
Then, for any  $h\in(0,1)$, $\delta$ sufficiently small, and $n \geq 1$ such that $nh^d \geq S\log(2(2+\delta)/\delta)$, and $nh^d  \geq S|\log(h)|$, we have with probability at least $1-\exp(-2\Delta^2S)$,
\begin{equation*}
    \|\hat{f}_{MoM} - f\|_\infty \leq C_1\sqrt{\frac{S \log\left(\frac{2(2+\delta)}{\delta}\right)|\log(h)|}{nh^d}} + C_2h^\alpha \;,
\end{equation*}
where $C_2=L\displaystyle\int\|u\|^\alpha K(u)du<\infty$ and $C_1$ is a constant that only depends on $\|f\|_\infty$, the dimension $d$, and the kernel properties.}

\begin{proof} 
From the definition of the MoM-KDE, we have the following implication \citep{lecue2018robust}
\begin{equation*}
     \left\{\sup_x \left\lvert \hat{f}_{MoM}(x) - f(x) \right\rvert \geq \varepsilon \right\} \Longrightarrow \left\{\sup_x \sum^S_{k=1} I\left(\left\lvert \hat{f}_{n_s}(x) - f(x) \right\rvert > \varepsilon \right) \geq S/2  \right\} \; .
\end{equation*}
Thus to upper-bound the probability of the left-hand event, it suffices to upper-bound the probability of the right-hand event.
Moreover, we have
\begin{align*}
    &\left\lvert \hat{f}_{n_s}(x) - f(x) \right\rvert \leq \sup_x \left\lvert \hat{f}_{n_s}(x) - f(x) \right\rvert \\
    \Longrightarrow \quad & I\left(\left\lvert \hat{f}_{n_s}(x) - f(x) \right\rvert > \varepsilon \right) \leq I\left(\sup_x \left\lvert \hat{f}_{n_s}(x) - f(x) \right\rvert > \varepsilon \right) \\
    \Longrightarrow \quad & \sum^S_{k=1} I\left(\left\lvert \hat{f}_{n_s}(x) - f(x) \right\rvert > \varepsilon \right) \leq \sum^S_{s=1} I\left(\sup_x \left\lvert \hat{f}_{n_s}(x) - f(x) \right\rvert > \varepsilon \right) \\
    \Longrightarrow \quad & \sup_x \sum^S_{s=1} I\left(\left\lvert \hat{f}_{n_s}(x) - f(x) \right\rvert > \varepsilon \right) \leq \sum^S_{s=1} I\left(\sup_x \left\lvert \hat{f}_{n_s}(x) - f(x) \right\rvert > \varepsilon \right) \; ,
\end{align*}
which implies that 
\begin{equation*}
\IP\left(\sup_x \sum^S_{s=1} I\left(\left\lvert \hat{f}_{n_s}(x) - f(x) \right\rvert > \varepsilon \right)\geq S/2\right) \leq \IP\left(\sum^S_{s=1} I\left(\sup_x \left\lvert \hat{f}_{n_s}(x) - f(x) \right\rvert > \varepsilon \right) \geq S/2\right) \; .
\end{equation*}

Let $Z_s=I\left(\sup_x \left\lvert \hat{f}_{n_s}(x) - f(x) \right\rvert > \varepsilon \right)$ and let $\calS = \big\{s \in \{1, \cdots, S\} \mid B_s \cap \calO = \emptyset  \big\}$ i.e. the set of indices $s$ such that the block $B_s$ does not contain any outliers. Since $\displaystyle\sum_{s \in \mathcal{S}^C} I(\cdot)$ is bounded by $\lvert \mathcal{O} \rvert$, almost surely, the following holds.

\begin{align}
    &\sum^S_{s=1} I\left(\sup_x \left\lvert \hat{f}_{n_s}(x) - f(x) \right\rvert > \varepsilon \right)= \sum^S_{s=1} Z_s
    = \sum_{s\in \calS} Z_s + \sum_{s\in \calS^C} Z_s \nonumber \\
    &\leq \sum_{s\in \calS} Z_s + \lvert \mathcal{O} \rvert \nonumber \\
    &=\sum_{s\in \calS} [ Z_k - \IE\left(Z_s\right) + \IE\left(Z_s\right)  ] + \lvert \mathcal{O} \rvert \nonumber \\
    &= \sum_{s\in \calS} [ Z_s - \IE\left(Z_s\right)] + \sum_{s\in \calS} \IE\left(Z_s\right)  + \lvert \mathcal{O} \rvert \nonumber \\
    &\leq \sum^S_{s= 1} [ Z_s - \IE\left(Z_s\right)] + S \cdot \IE\left(Z_1\right)  + \lvert \mathcal{O} \rvert\nonumber \\
    \label{eq:last_line_prop}
    &\leq \sum^S_{s= 1} [ Z_s - \IE\left(Z_s\right)] + S \cdot \IP\left(\sup_x \left\lvert \hat{f}_{n_1}(x) - f(x) \right\rvert > \varepsilon\right)  + \lvert \mathcal{O} \rvert \; ,
\end{align}
where $Z_1$ is assumed, without loss of generality, to be associated to a block not containing outliers. This block always exists thanks to the hypothesis $S>(2+\delta)|\mathcal{O}|$. \\

Let $\varepsilon = C_1\sqrt{\frac{S\log(\frac{2(2+\delta)}{\delta})|\log(h)|}{nh^d}} + C_2h^\alpha$, then using Lemma $1$ with $S = \log(\frac{2(2+\delta)}{\delta})$, we have

\begin{equation*}
    \IP\left(\sup_x \left\lvert \hat{f}_{n_1}(x) - f(x) \right\rvert > \varepsilon\right) \leq \frac{\delta}{2(2+\delta)} \; .
\end{equation*}

Combining this last inequality with equation \eqref{eq:last_line_prop} leads to

\begin{align*}
    \IP\left(\sum^S_{s=1} I\left(\sup_x \left\lvert \hat{f}_{n_s}(x) - f(x) \right\rvert > \varepsilon \right) \geq S/2\right)  & \leq \IP\left(\sum^S_{s= 1} [ Z_s - \IE\left(Z_s\right)] + S \cdot \frac{\delta}{2(2+\delta}  + \lvert \mathcal{O}\rvert \geq S/2\right) \\
    & \leq \IP\left(\sum^S_{s= 1} [ Z_s - \IE\left(Z_s\right)] \geq S\left(\frac{1}{2} - \frac{\delta}{2(2+\delta)} - \frac{|\calO|}{S}\right)\right) \\
    & \leq \IP\left(\sum^S_{s= 1} [ Z_s - \IE\left(Z_s\right)] \geq S\left( \frac{1}{2+\delta} - \frac{|\calO|}{S}\right)\right)
\end{align*}

Tacking $\Delta = \left(\frac{1}{2+\delta} - \frac{|\calO|}{S}\right) > 0$ and applying Hoeffding’s inequality to the right-hand side of the previous equation gives

\begin{align*}
    \IP\left(\sum^S_{s=1} I\left(\sup_x \left\lvert \hat{f}_{n_s}(x) - f(x) \right\rvert > \varepsilon \right) \geq S/2 \right) & \leq e^{-2S \Delta^2} \; ,
\end{align*}
which concludes the proof.
\end{proof}

\paragraph{Proposition 2.}\emph{($L_1$-consistency in probability) If $n/S\rightarrow \infty$, $h\rightarrow 0$, $n h^d \rightarrow \infty$, and $S>2|\mathcal{O}|,$ then
\begin{equation*}
  \|\hat{f}_{MoM} - f\|_1 \overset{\mathcal{P}}{\underset{n\rightarrow \infty}{\longrightarrow}} 0 \; .
\end{equation*}}

\begin{proof}
We first rewrite the MoM-KDE as
\begin{equation*}
    \hat{f}_{MoM}(x) = \sum^S_{s=1} \hat{f}_{n_s}(x) I_{A_s}(x) \; ,
\end{equation*}
where $A_s = \left\{x \mid \hat{f}_{MoM}(x) = \hat{f}_{n_s}(x) \right\}$. Without loss of generality, we assume that
\begin{equation*}
    A_k \overset{S}{\underset{s \neq \ell}{\cap}} A_\ell = \emptyset,  \quad \overset{S}{\underset{s = 1}{\cup}} A_s = \IR^d, \quad \text{ and } \quad \sum^S_{s=1} I_{A_s}(x) = 1 \; .
\end{equation*}

\begin{align}
    \int \left\lvert \hat{f}_{MoM}(x) - f(x) \right\rvert dx & = \int \left\lvert \sum^S_{s=1} \hat{f}_{n_s}(x) I_{A_s}(x) - f(x) \right\rvert dx \nonumber\\
    &= \int \left\lvert \sum^S_{s=1} \left(\hat{f}_{n_s}(x) - f(x)\right)I_{A_s}(x) \right\rvert dx \nonumber\\
    &\leq \int \sum^S_{s=1} \left\lvert \hat{f}_{n_s}(x) - f(x) \right\rvert I_{A_s}(x) dx  \nonumber\\
    &= \sum^S_{s=1} \int_{A_s} \left\lvert \hat{f}_{n_s}(x) - f(x) \right\rvert dx  \nonumber\\
    \label{eq:prop2}
    &= \sum_{s \in \calS} \int_{A_s} \left\lvert \hat{f}_{n_s}(x) - f(x) \right\rvert dx + \sum_{s \in \calS^C} \int_{A_s} \left\lvert \hat{f}_{n_s}(x) - f(x) \right\rvert dx  \; .
\end{align}
From the $L_1$-consistency of the KDE in probability, if the number of anomalies grows at a small enough speed \citep{devroyegy}, the left part is bounded, i.e.
\begin{equation}
    \sum_{s \in \calS} \int_{A_s} \left\lvert \hat{f}_{n_s}(x) - f(x) \right\rvert dx \leq \sum_{s \in \calS} \int \left\lvert \hat{f}_{n_s}(x) - f(x) \right\rvert dx \overset{\mathcal{P}}{\underset{n \rightarrow \infty}{\longrightarrow}} 0 \; .
    \label{eq:right_const}
\end{equation}

We now upper-bound the right part of equation \eqref{eq:prop2}. Let consider a particular block $A_s$ where $s \in \calS^C$. In this block, the estimator $f_{n_s}$ is selected and is calculated with samples containing anomalies. As $\forall x \in A_s$, $f_{n_s}(x)$ is the median (by definition), if $S > 2 \lvert \calO \rvert$, we can always find a $s^\prime\in \calS$ such that $f_{n_s}(x)\leq f_{n_{s^\prime}}(x)$ or $f_{n_s}(x)\geq f_{n_{s^\prime}}(x)$. \\

\noindent
Now let denote by ${A^+_s} = \left\{x \in A_s \mid \hat{f}_{n_s}(x) \geq f(x)  \right\}$ and ${A^-_s} = \left\{x \in A_s \mid \hat{f}_{n_s}(x) < f(x)  \right\}$. We have ${A^+_s}\cup{A^-_s} = A_s$ and each one of these blocks can be decomposed respectively into $|\calS|$ sub-blocks (not necessarily disjoint) $\{A^{s^\prime,+}_s\}_{s^\prime \in \calS}$ and $\{A^{s^\prime,-}_s\}_{s^\prime \in \calS}$ such that $\forall s^\prime \in \calS$,
\small
${A^{s^\prime,+}_s} = \left\{x \in A_s \mid \hat{f}_{n_{s^\prime}}(x) \geq \hat{f}_{n_s}(x) \geq f(x)  \right\}$ \normalsize and \small ${A^{s^\prime,-}_s} = \left\{x \in A_s \mid \hat{f}_{n_{s^\prime}}(x) \leq \hat{f}_{n_s}(x) < f(x)  \right\}$.
\normalsize
Finally, the right-hand term of equation \eqref{eq:prop2} can be upper-bounded by

\begin{align*}
     \sum_{s \in \calS^C} \int_{A_s} \left\lvert \hat{f}_{n_s}(x) - f(x) \right\rvert dx & \leq  \sum_{s \in \calS^C} \int_{A^+_s} \left\lvert \hat{f}_{n_s}(x) - f(x) \right\rvert dx + \int_{A^-_s} \left\lvert \hat{f}_{n_s}(x) - f(x) \right\rvert dx \\
     & \leq  \sum_{s \in \calS^C}\sum_{s^\prime \in \calS} \int_{A^{s^\prime,+}_s} \left\lvert \hat{f}_{n_s}(x) - f(x) \right\rvert dx + \int_{A^{s^\prime,-}_s} \left\lvert \hat{f}_{n_s}(x) - f(x) \right\rvert dx \\
     & \leq  \sum_{s \in \calS^C}\sum_{s^\prime \in \calS} \int_{A^{s^\prime,+}_s} \left\lvert \hat{f}_{n_{s^\prime}}(x) - f(x) \right\rvert dx + \int_{A^{s^\prime,-}_s} \left\lvert \hat{f}_{n_{s^\prime}}(x) - f(x) \right\rvert dx \\
     & \leq  \sum_{s \in \calS^C}\sum_{s^\prime \in \calS} \int \left\lvert \hat{f}_{n_{s^\prime}}(x) - f(x) \right\rvert dx + \int \left\lvert \hat{f}_{n_{s^\prime}}(x) - f(x) \right\rvert dx  \; .
\end{align*}

Since $\forall s^\prime \in \calS$ we have $\displaystyle\int \lvert \hat{f}_{n_{s^\prime}}(x) - f(x) \rvert dx \overset{\mathcal{P}}{\underset{n \rightarrow \infty}{\longrightarrow}} 0$, we can conclude using similar arguments as those used for \eqref{eq:right_const} that $\displaystyle\sum_{s \in \calS^C} \int_{A_s} \lvert \hat{f}_{n_s}(x) - f(x) \rvert dx\overset{\mathcal{P}}{\underset{n \rightarrow \infty}{\longrightarrow}} 0$, which concludes the proof.

\end{proof}

\section{Additional results}

\begin{figure}
	\centering
	$D_{\text{KL}}(f, \hat{f})$\\
	\medskip
	\begin{minipage}[t]{\ratio\linewidth}
		\centering
		\includegraphics[height=\ratiob\linewidth]{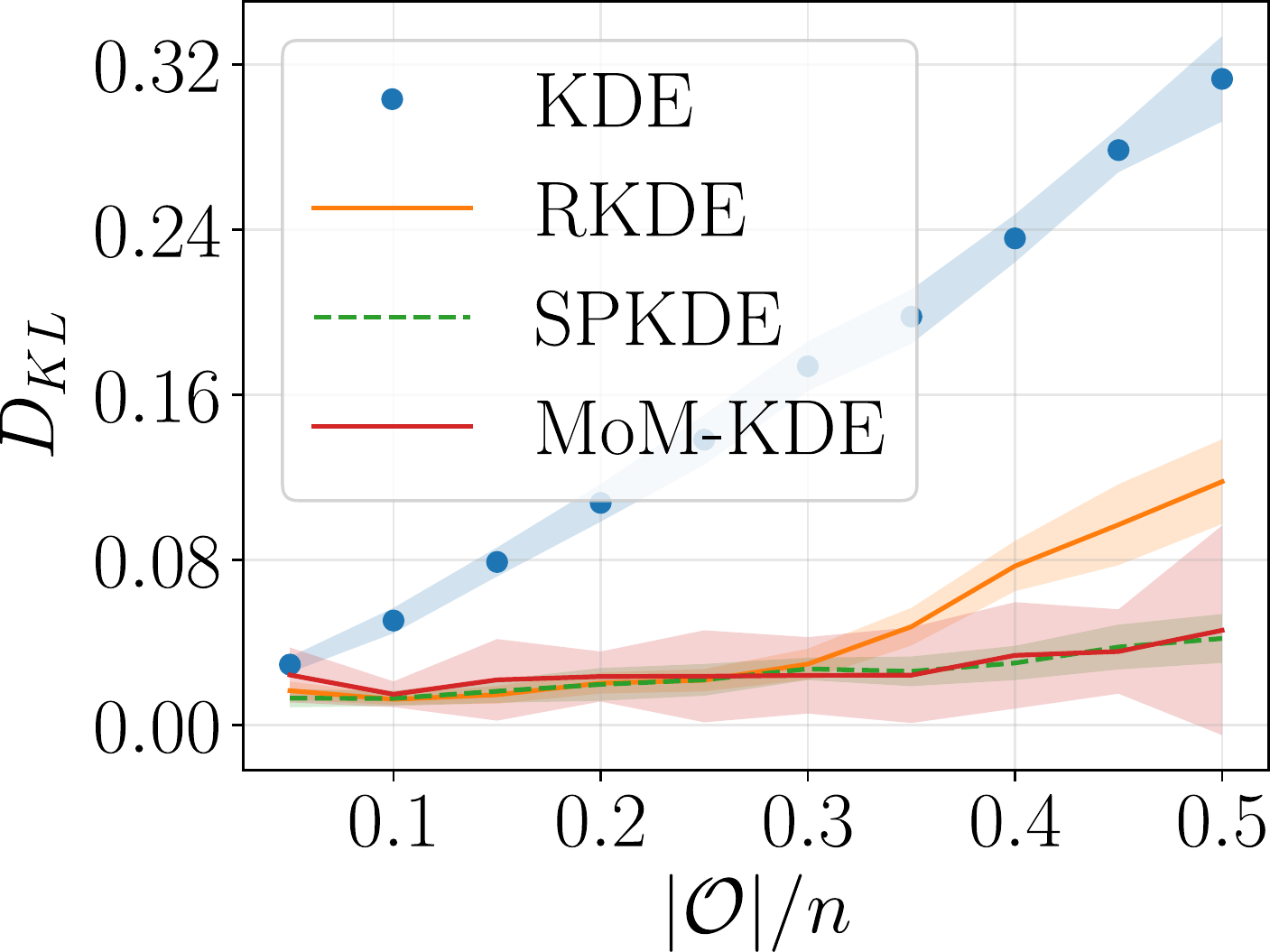}
	\end{minipage}\hfill
	\begin{minipage}[t]{\ratio\linewidth}
		\centering
		\includegraphics[height=\ratiob\linewidth]{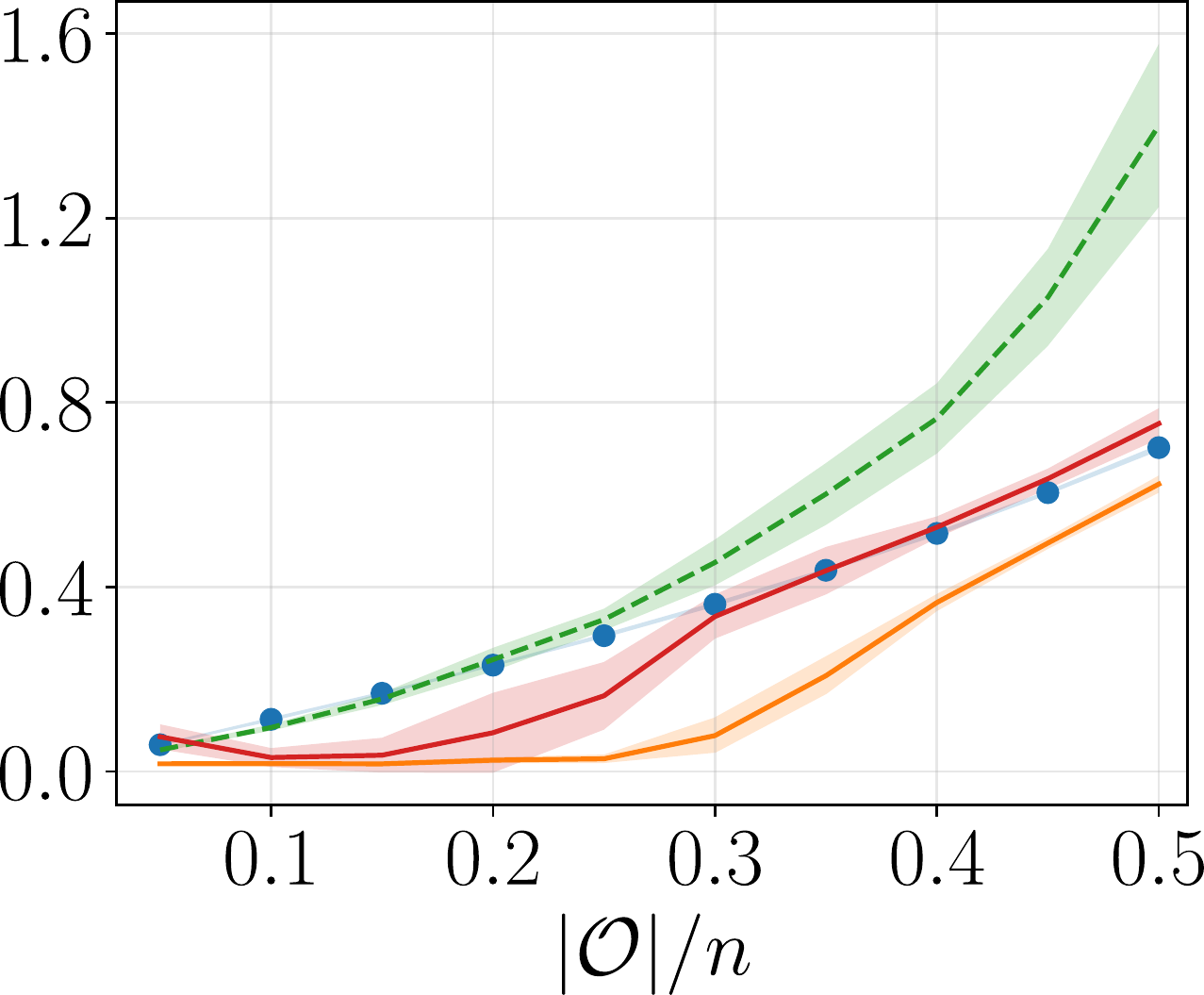}
	\end{minipage}
	\begin{minipage}[t]{\ratio\linewidth}
		\centering
		\includegraphics[height=\ratiob\linewidth]{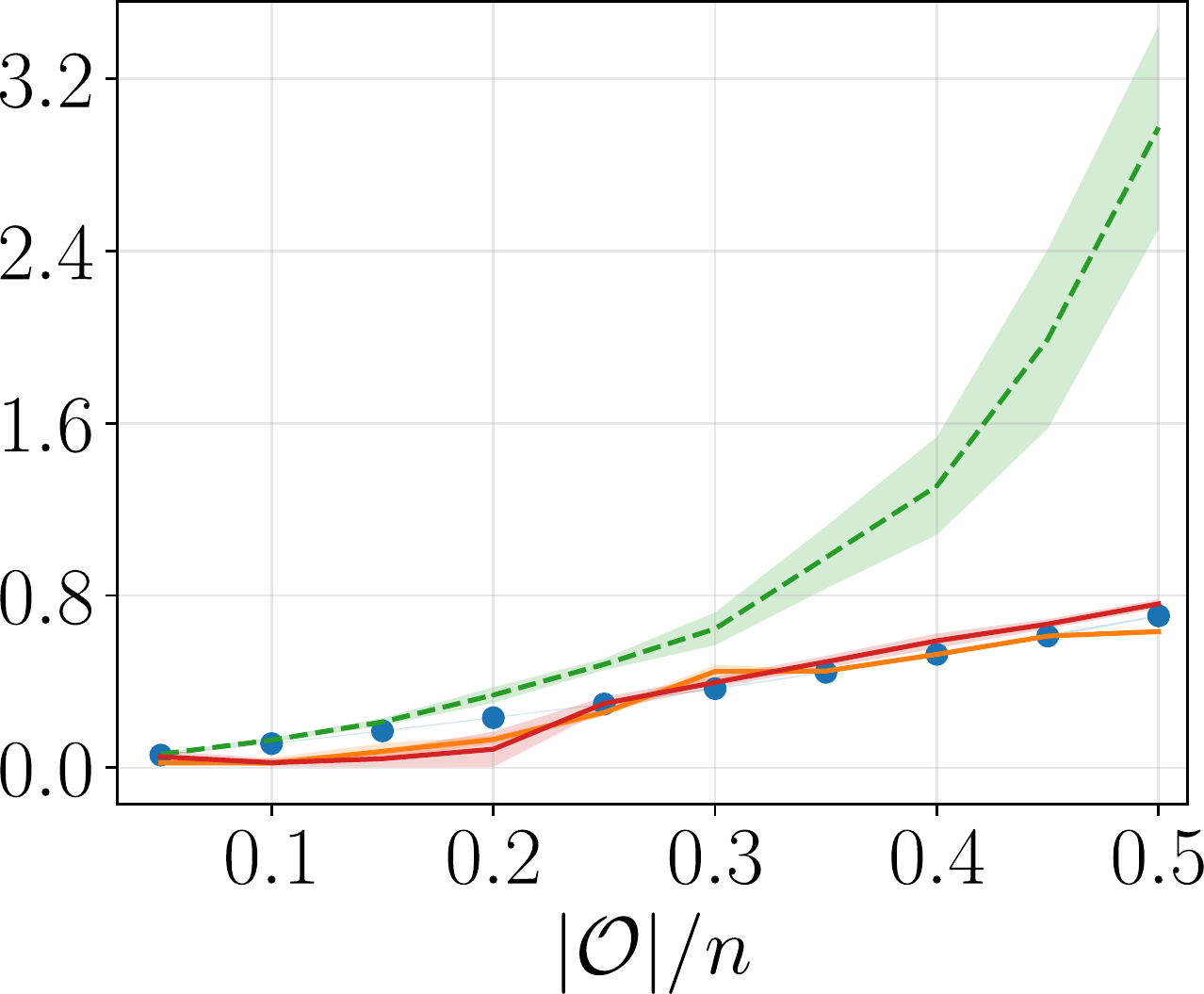}
	\end{minipage}
	\begin{minipage}[t]{\ratio\linewidth}
		\centering
		\includegraphics[height=\ratiob\linewidth]{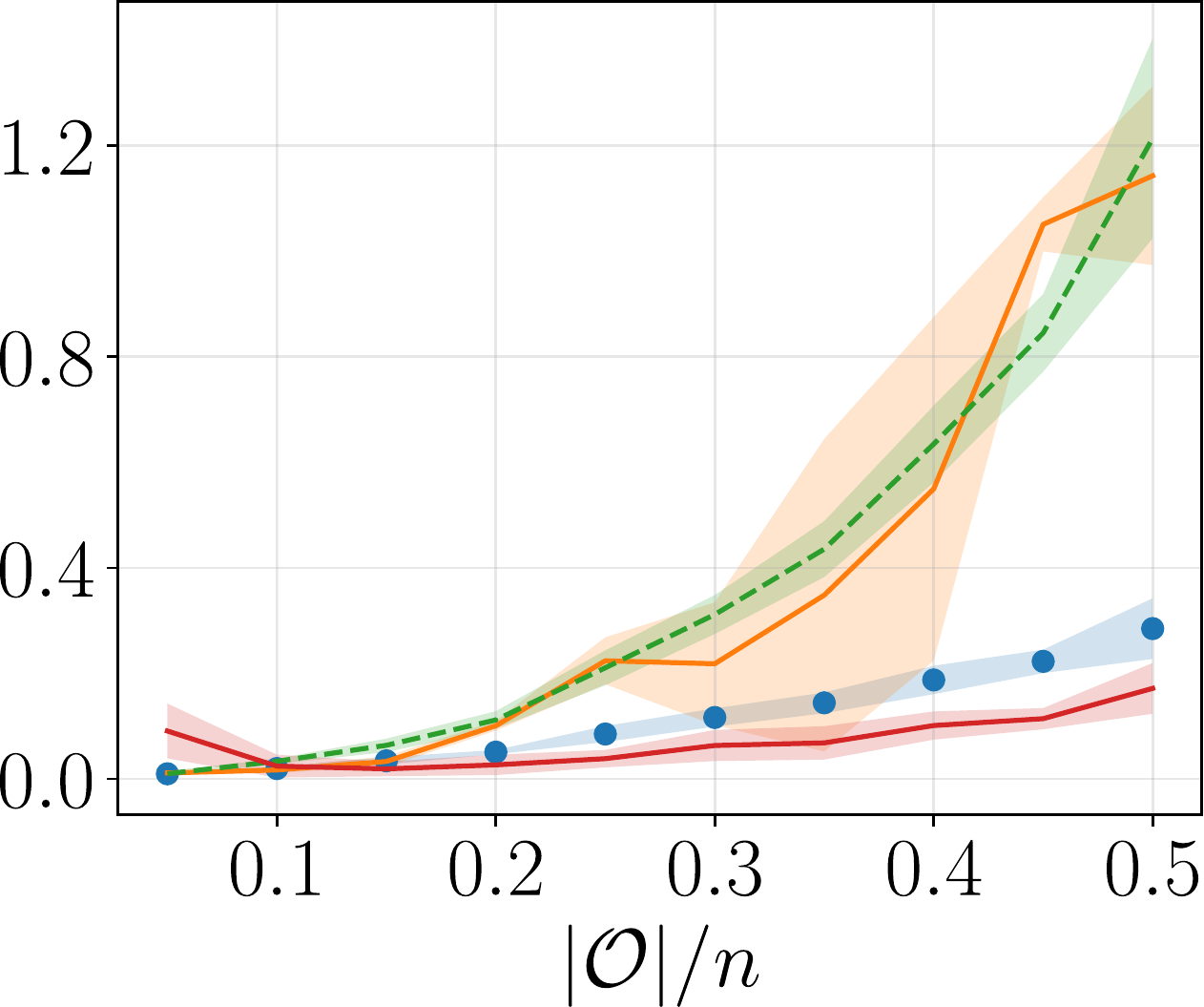}
	\end{minipage}\\
	$D_{\text{KL}}(\hat{f}, f)$\\
	\medskip
	\renewcommand{\ratiob}{0.82}
	\begin{minipage}[t]{\ratio\linewidth}
		\centering
		\includegraphics[height=0.83\linewidth]{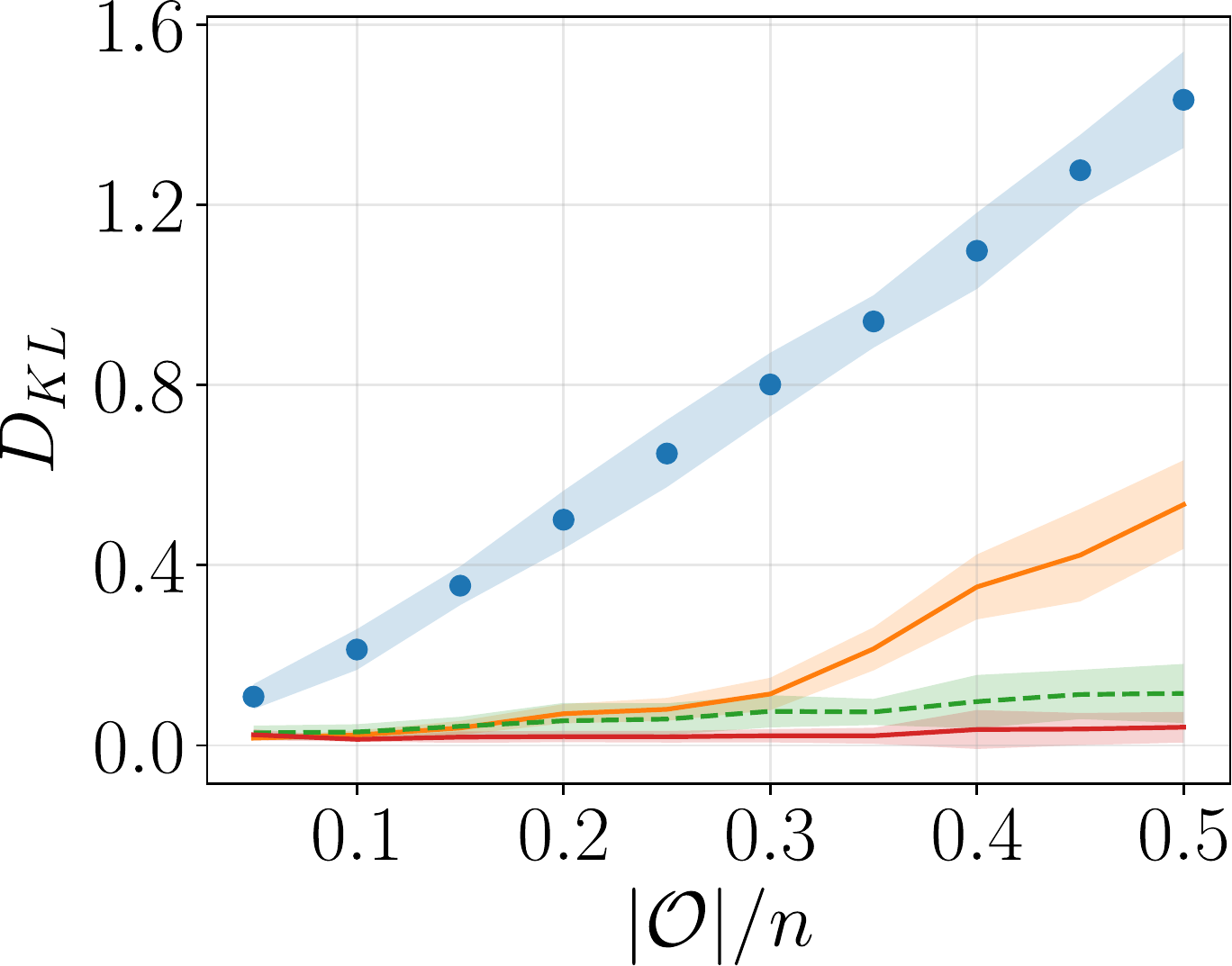}
	\end{minipage}\hfill
	\begin{minipage}[t]{\ratio\linewidth}
		\centering
		\includegraphics[height=\ratiob\linewidth]{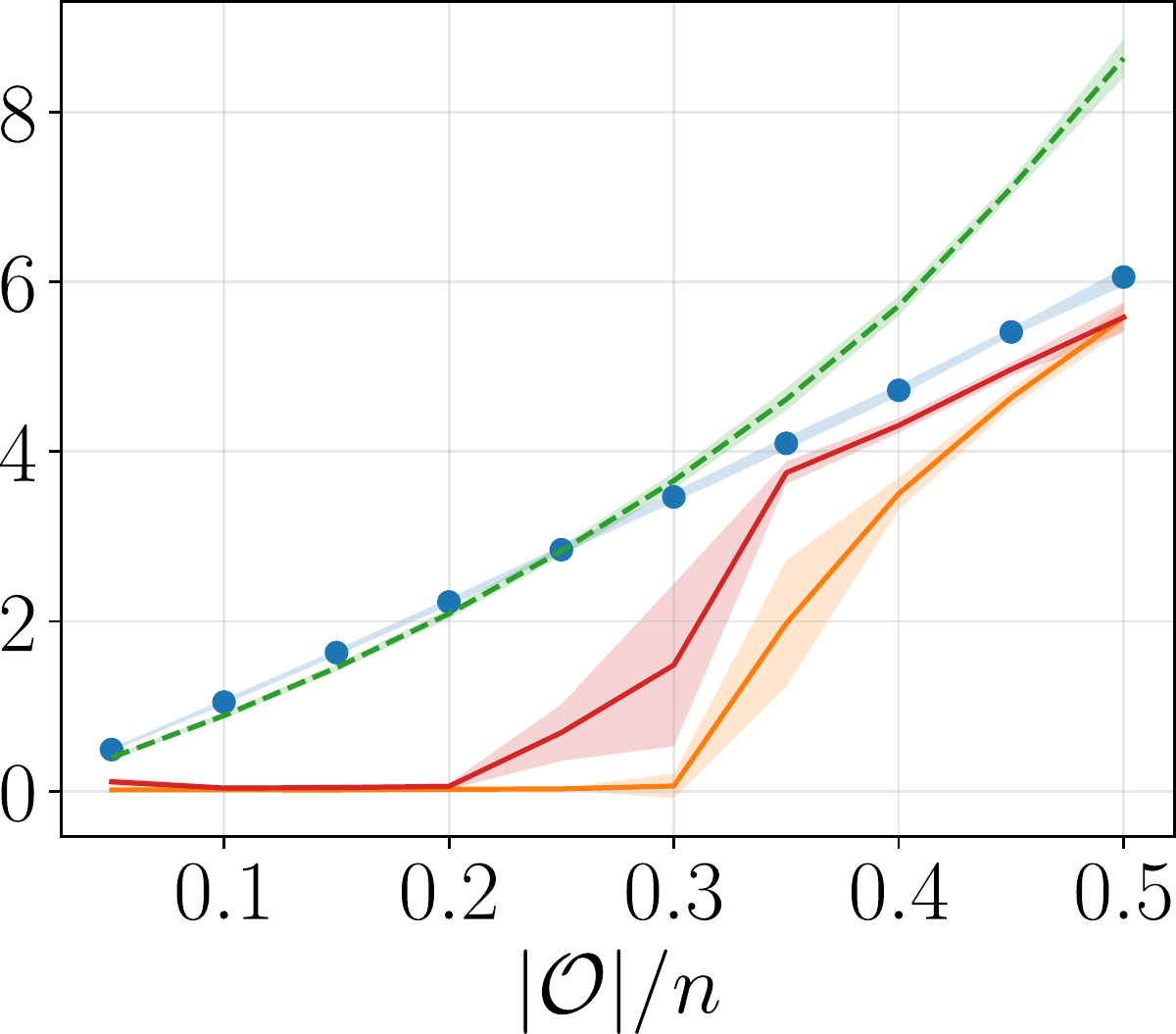}
	\end{minipage}
	\begin{minipage}[t]{\ratio\linewidth}
		\centering
		\includegraphics[height=\ratiob\linewidth]{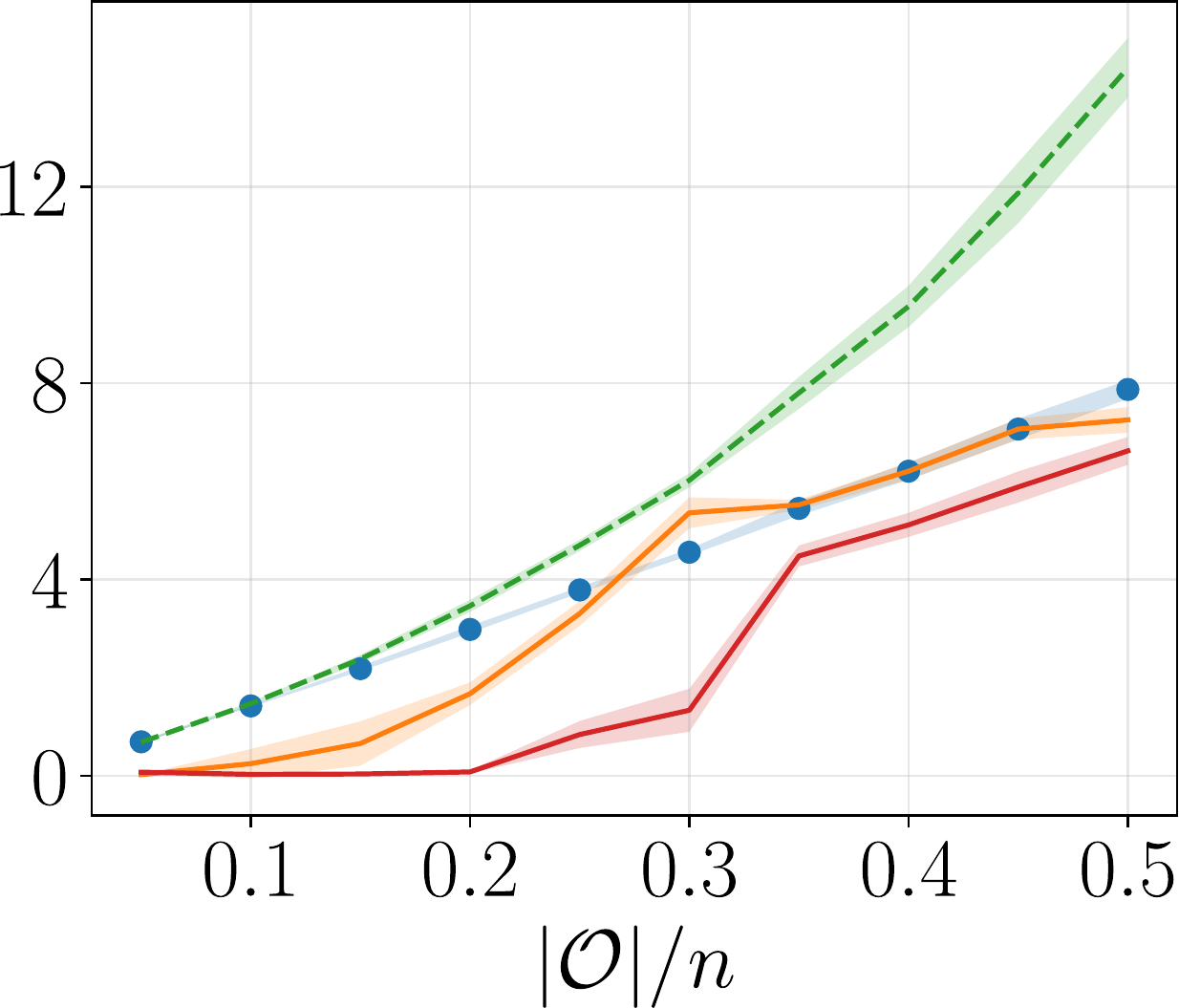}
	\end{minipage}
	\begin{minipage}[t]{\ratio\linewidth}
		\centering
		\includegraphics[height=\ratiob\linewidth]{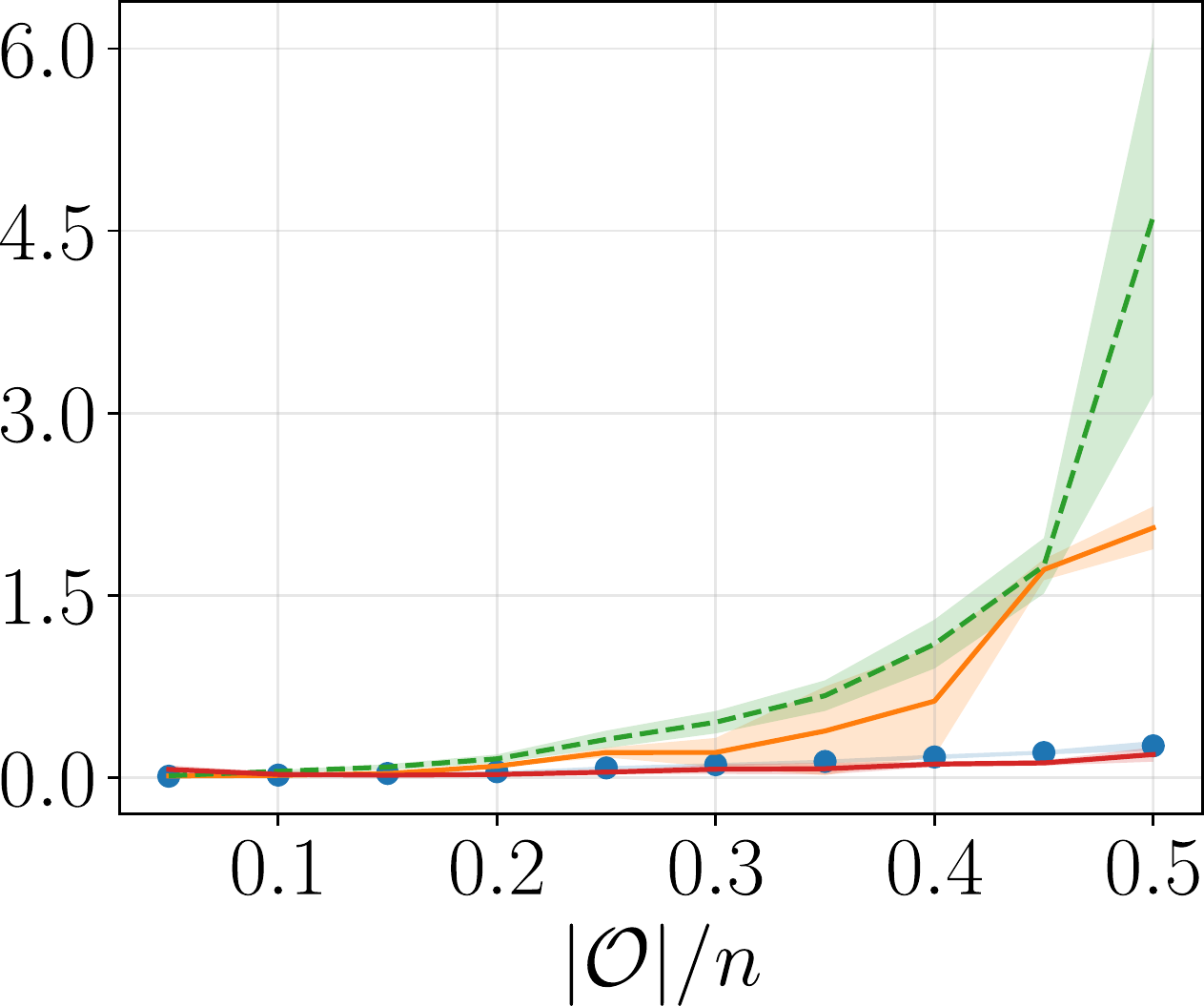}
	\end{minipage}
	\renewcommand{\ratiob}{0.8}
	AUC\\
	\medskip
	\begin{minipage}[t]{\ratio\linewidth}
		\centering
		\includegraphics[height=\ratiob\linewidth]{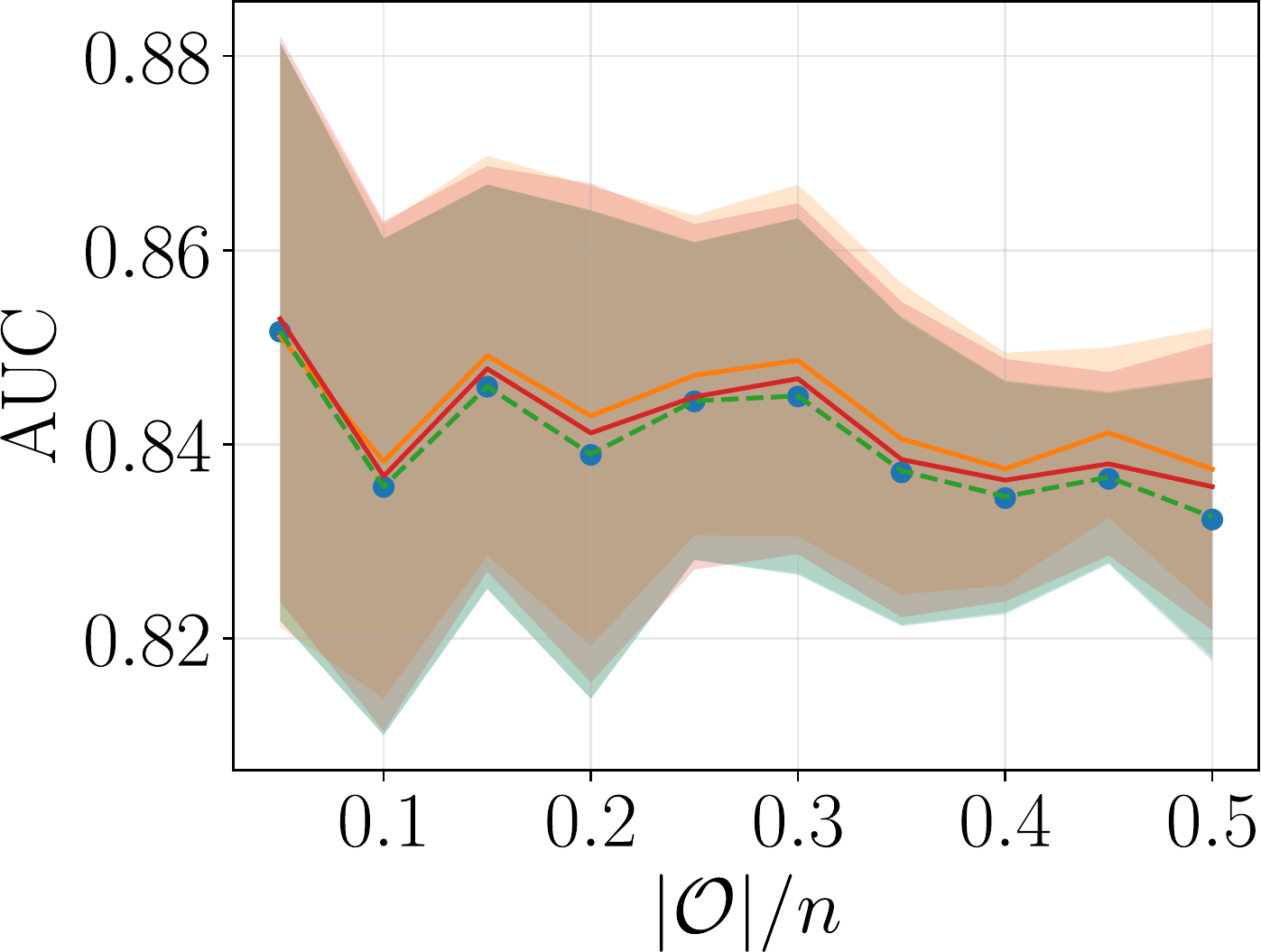}\\[2pt]
		\hspace*{25pt}{\small (a)\; Uniform}
	\end{minipage}\hfill
	\begin{minipage}[t]{\ratio\linewidth}
		\centering
		\includegraphics[height=\ratiob\linewidth]{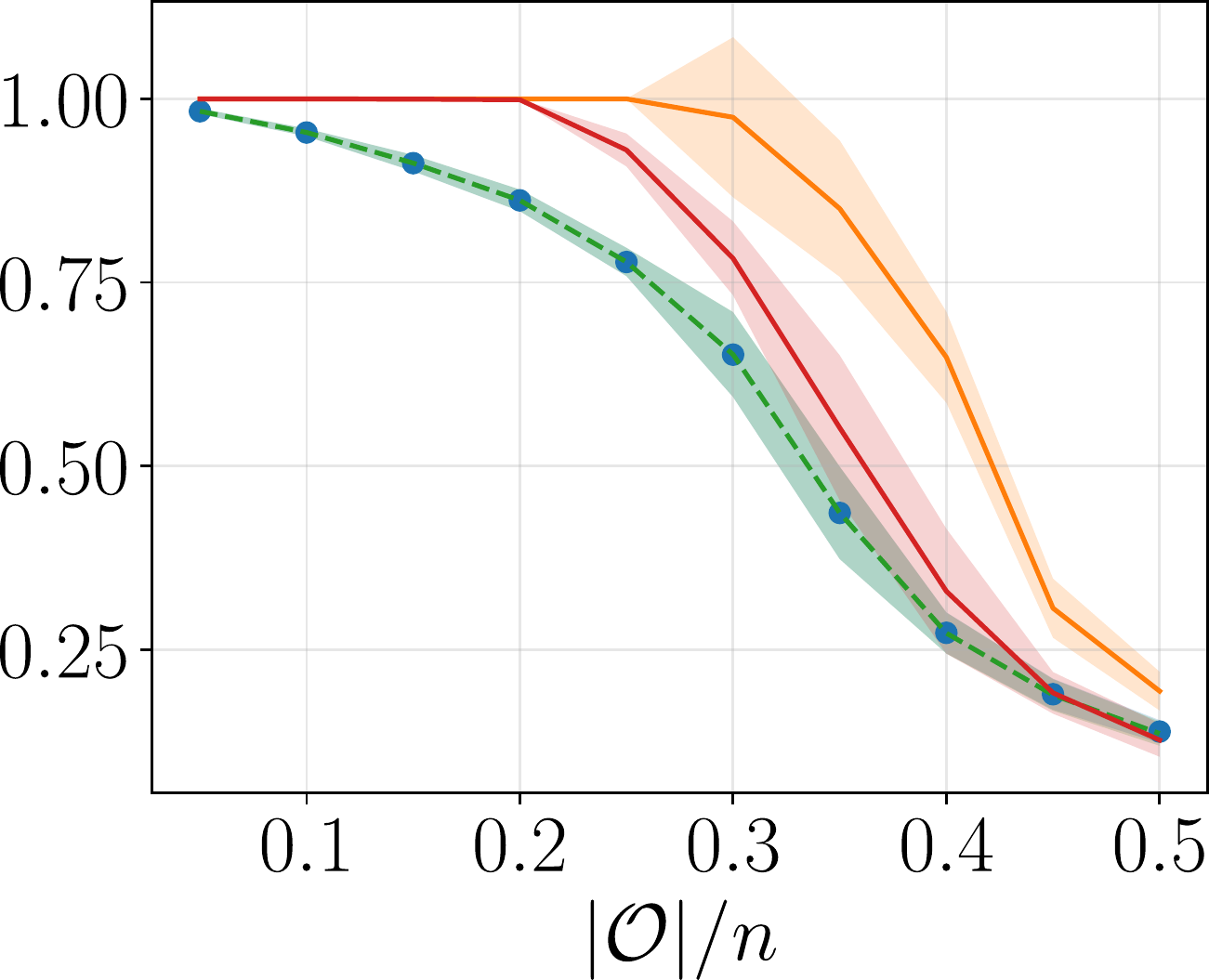}\\[2pt]
		{\small (b)\; Regular Gaussian}
	\end{minipage}
	\begin{minipage}[t]{\ratio\linewidth}
		\centering
		\includegraphics[height=\ratiob\linewidth]{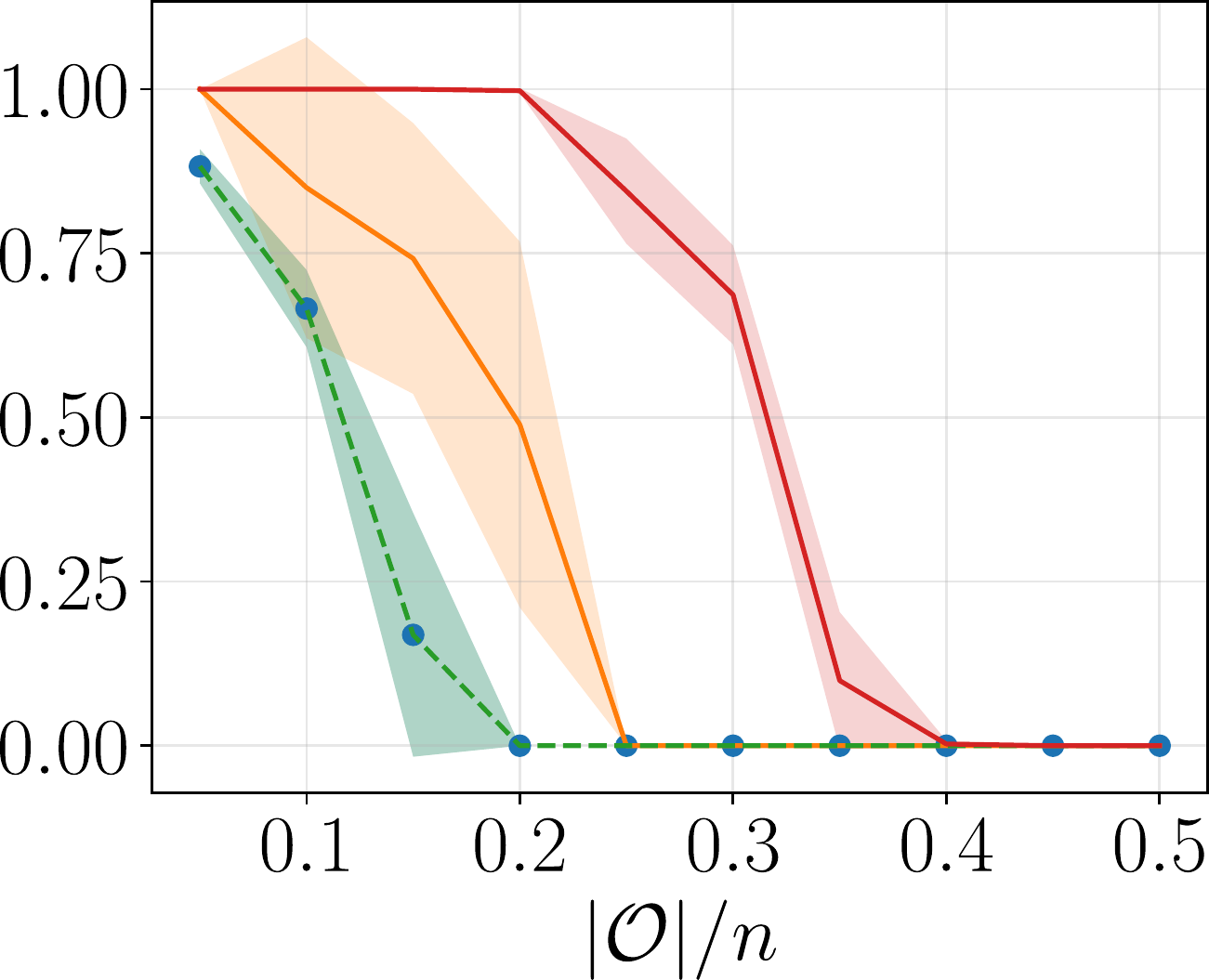}\\[2pt]
		{\small (c)\; Thin Gaussian}
	\end{minipage}
	\begin{minipage}[t]{\ratio\linewidth}
		\centering
		\includegraphics[height=\ratiob\linewidth]{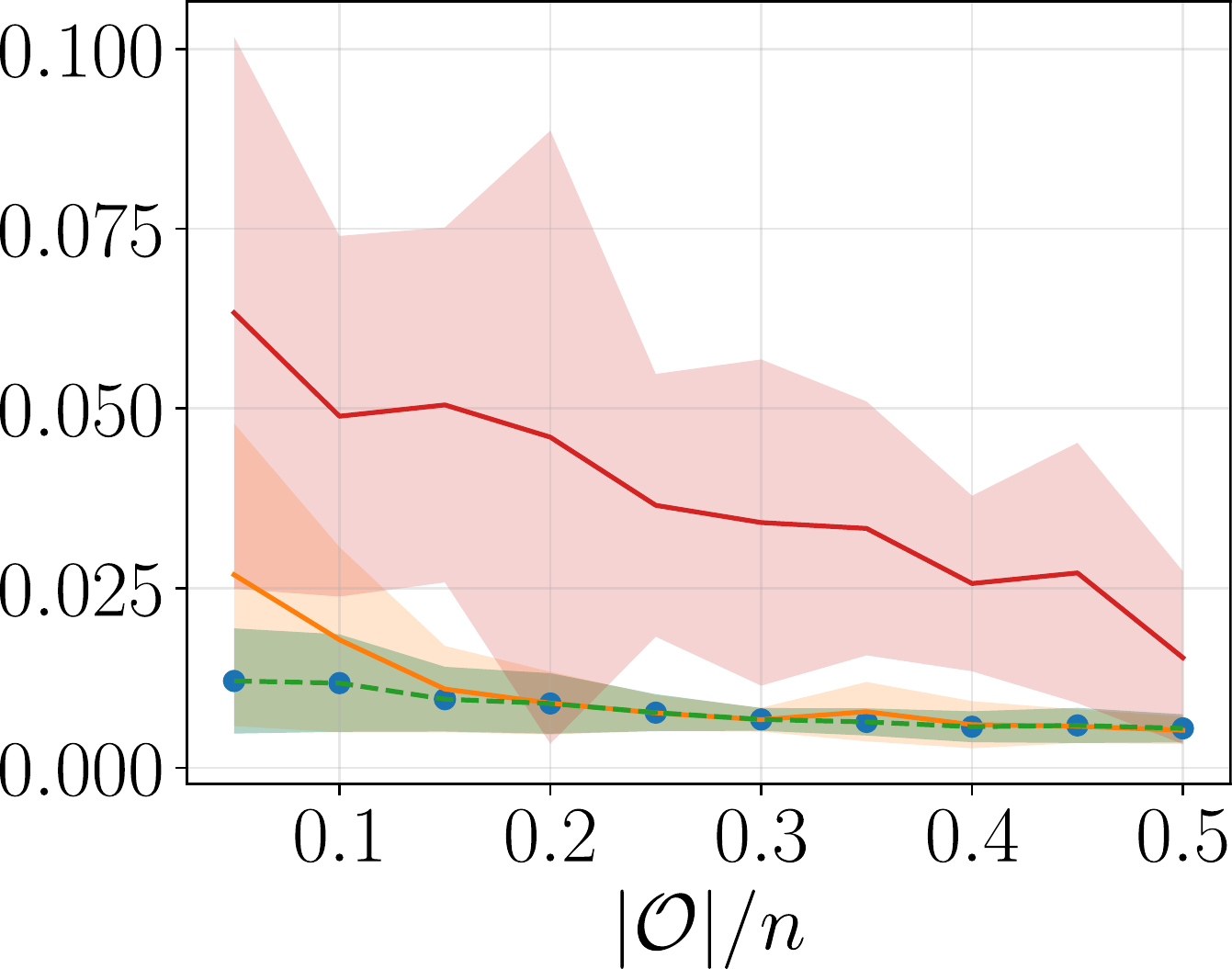}\\[2pt]
		{\small (d)\; Advsersarial Thin Gaussian}
	\end{minipage}
	\caption{Density estimation with synthetic data. The displayed metrics are the Kullback-Leibler divergence (a lower score means a better estimation of the true density) and the AUC (a higher score means a better detection of the outliers).}
	\label{fig:experiments_synth_app}
\end{figure}

As stated in the main paper, we display here the additional results containing:
\begin{itemize}[leftmargin=1cm]
	\item For synthetic data, the Kullback-Leibler divergence in both directions, i.e.\ $D_{\text{KL}}(\hat{f}, f)$ and $D_{\text{KL}}(f, \hat{f})$, and the ROC AUC measuring the performance of an anomaly detector based on $\hat{f}$. Results are displayed on Figure \ref{fig:experiments_synth_app}.
	\item For Digits data, the ROC AUC measuring the performance of an anomaly detector based on $\hat{f}$. As stated in the main paper, this is done under multiple scenarios, where outliers and inliers can be chosen among the nine available classes. Here we show the AUC when the outliers are set as one class (class 2 to class 9), and inliers are set as ``the rest'' of all classes. Results are displayed on Figure \ref{fig:experiments_real_app}.
\end{itemize}

\paragraph{Synthetic data.}
When considering the Kullback-Leibler divergence, results lead to a very similar conclusion as previously stated, that is, an overall good performance of MoM-KDE while its competitors, notably SPKDE, are more data-dependent.
When the density estimate $\hat{f}$ is used in a simple anomaly detector, results are quite different.
Indeed, when outliers are uniformly distributed, even if MoM-KDE seems to better estimate the true density (according to $D_{\text{JS}}$ and $D_{\text{KL}}$), this doesn't make $\hat{f}_{MoM}$ a better anomaly detector.
It seems that in this case, the outliers are either easily detected because distant from the density estimate, or located in dense regions, thus making them impossible to identify, and this for all density estimates provided by competitors.
In the case of adversarial contamination, the conclusion is quite similar.
Although MoM-KDE better fits the true density, the situation is extremely difficult for anomaly detection, hence making all competitors yield very poor results.
In the two other cases -- Gaussian outlier, anomaly detection results follow the density estimation.

\paragraph{Digits data.}
Results over Digits scenarios are inline with main conclusions over real data.
Although from one scenario to another, all methods have varied results, the overall observation is that MoM-KDE is either similar or better than its competitors.

\renewcommand{\ratio}{0.24}
\renewcommand{\ratiob}{0.83}
\begin{figure}
	\centering
	\begin{minipage}[t]{\ratio\linewidth}
		\centering
		\includegraphics[height=\ratiob\linewidth]{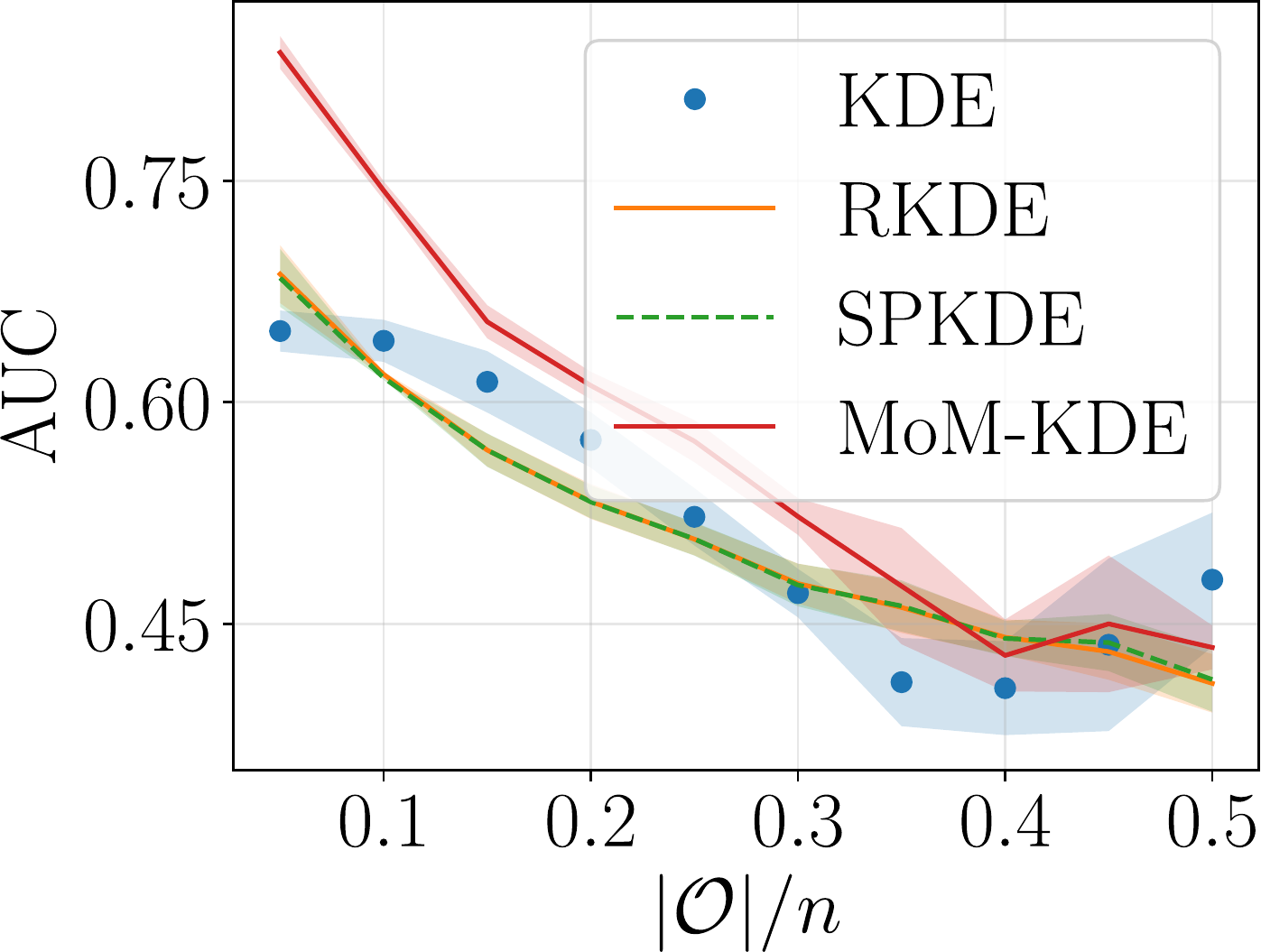}
		\hspace*{25pt}{\small(a)\; $\calO$: 2, $\calI$: all}
	\end{minipage}\hfill
	\begin{minipage}[t]{\ratio\linewidth}
		\centering
		\includegraphics[height=\ratiob\linewidth]{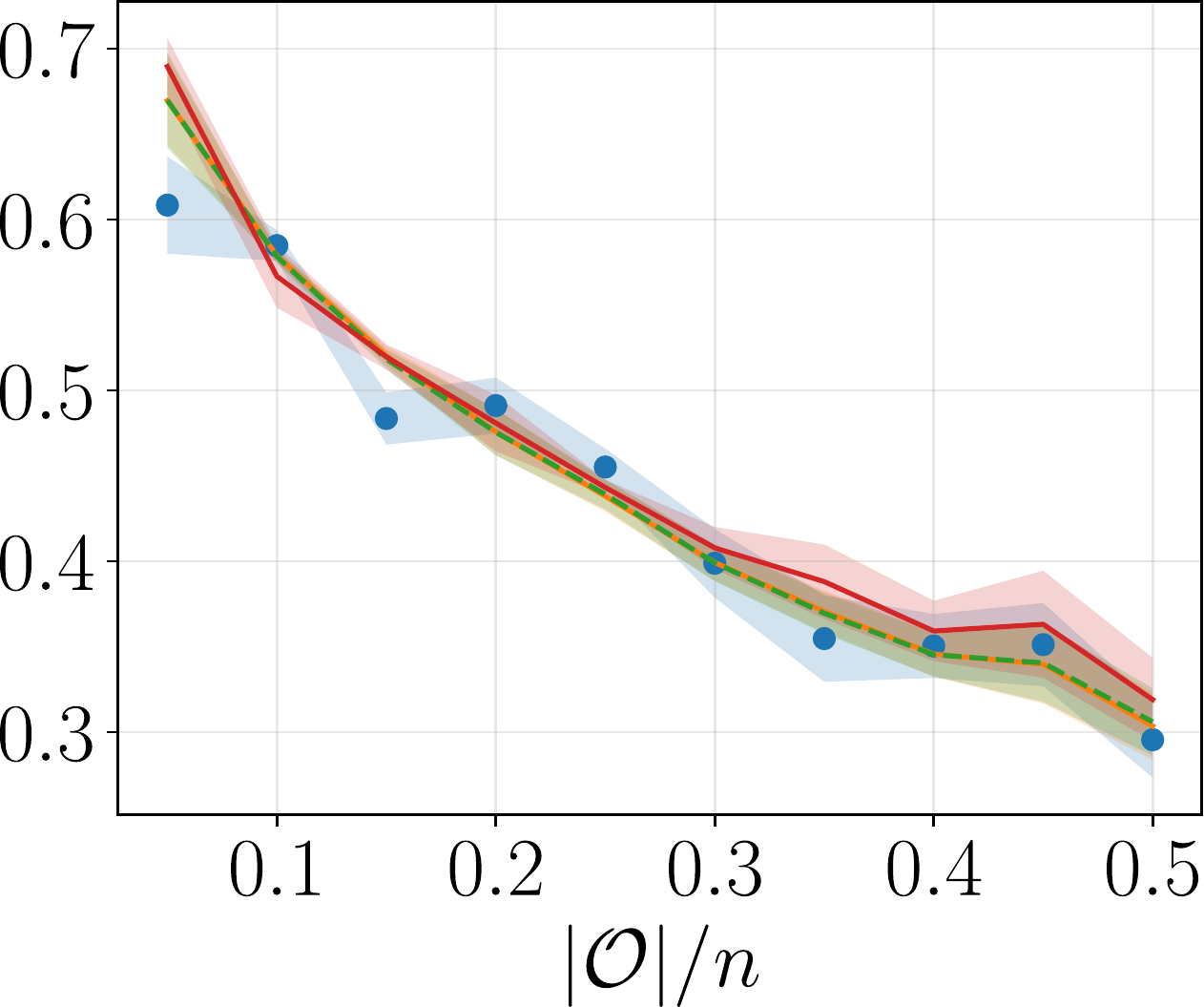}
		{\small (b)\; $\calO$: 3, $\calI$: all}
	\end{minipage}
	\begin{minipage}[t]{\ratio\linewidth}
		\centering
		\includegraphics[height=\ratiob\linewidth]{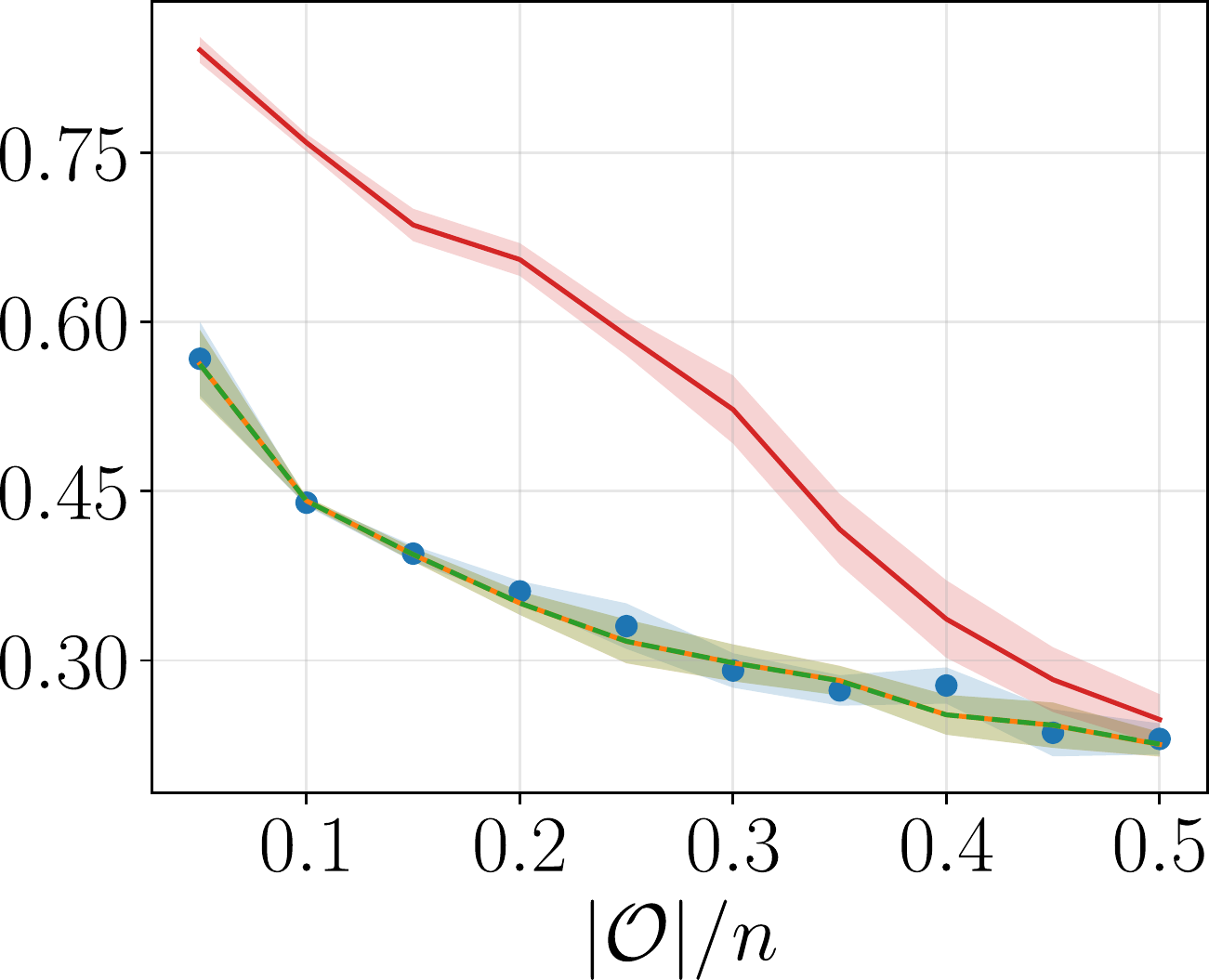}
		{\small (c)\; $\calO$: 4, $\calI$: all}
	\end{minipage}
	\begin{minipage}[t]{\ratio\linewidth}
		\centering
		\includegraphics[height=\ratiob\linewidth]{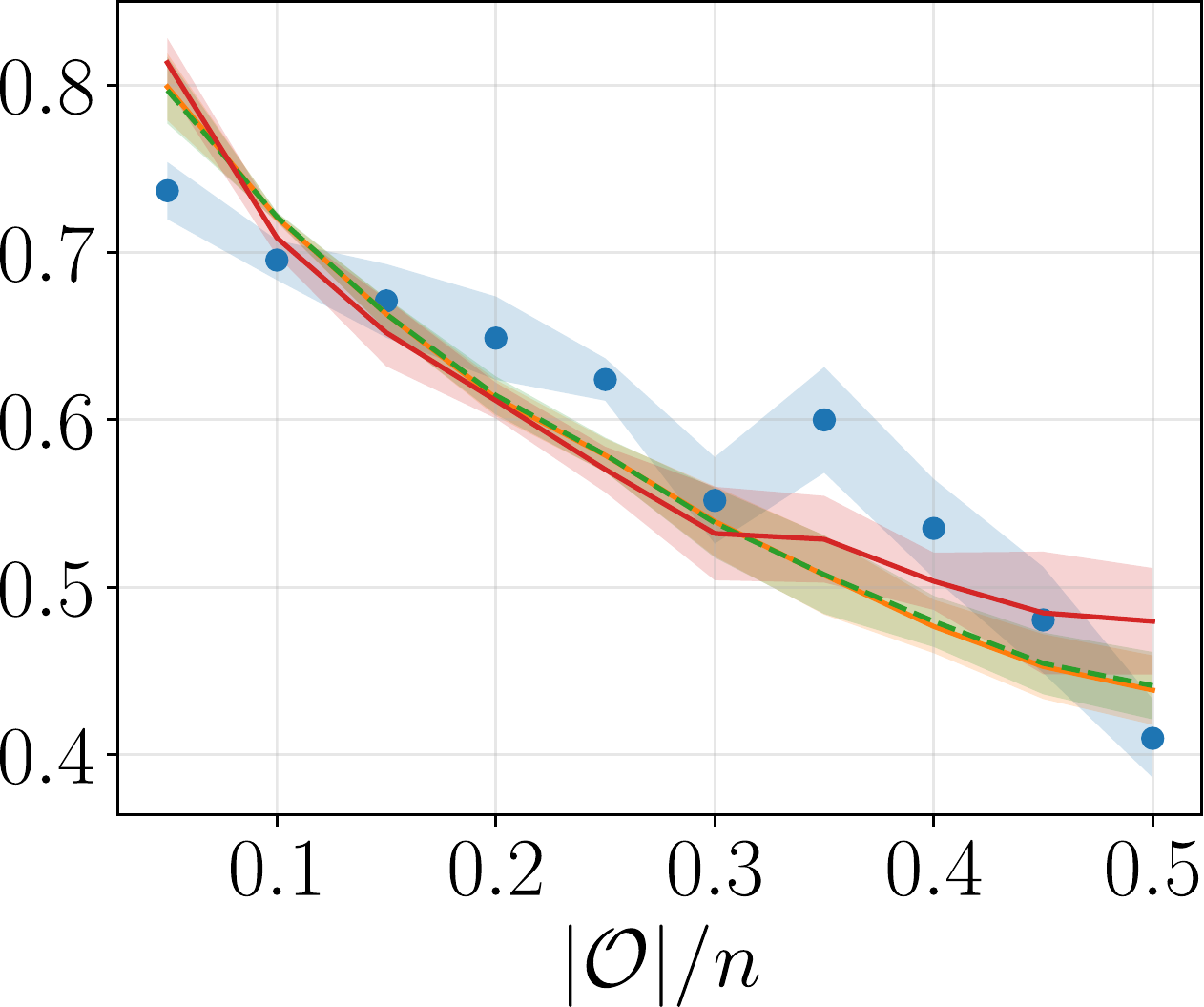}
		{\small (d)\; $\calO$: 5, $\calI$: all}
	\end{minipage}\\
	\medskip
	\begin{minipage}[t]{\ratio\linewidth}
		\centering
		\includegraphics[height=\ratiob\linewidth]{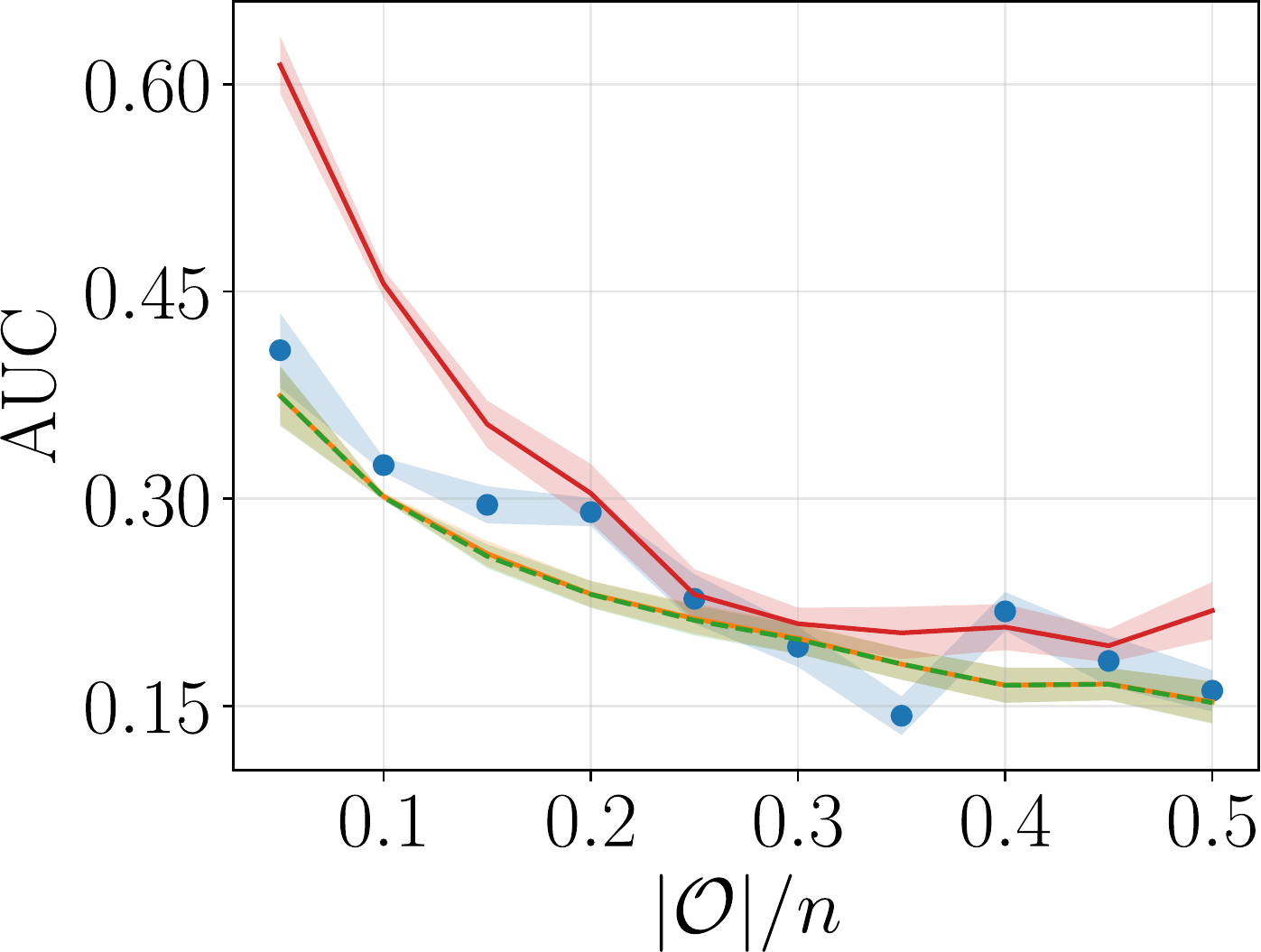}
		\hspace*{25pt}{\small (e)\; $\calO$: 6, $\calI$: all}
	\end{minipage}\hfill
	\begin{minipage}[t]{\ratio\linewidth}
		\centering
		\includegraphics[height=\ratiob\linewidth]{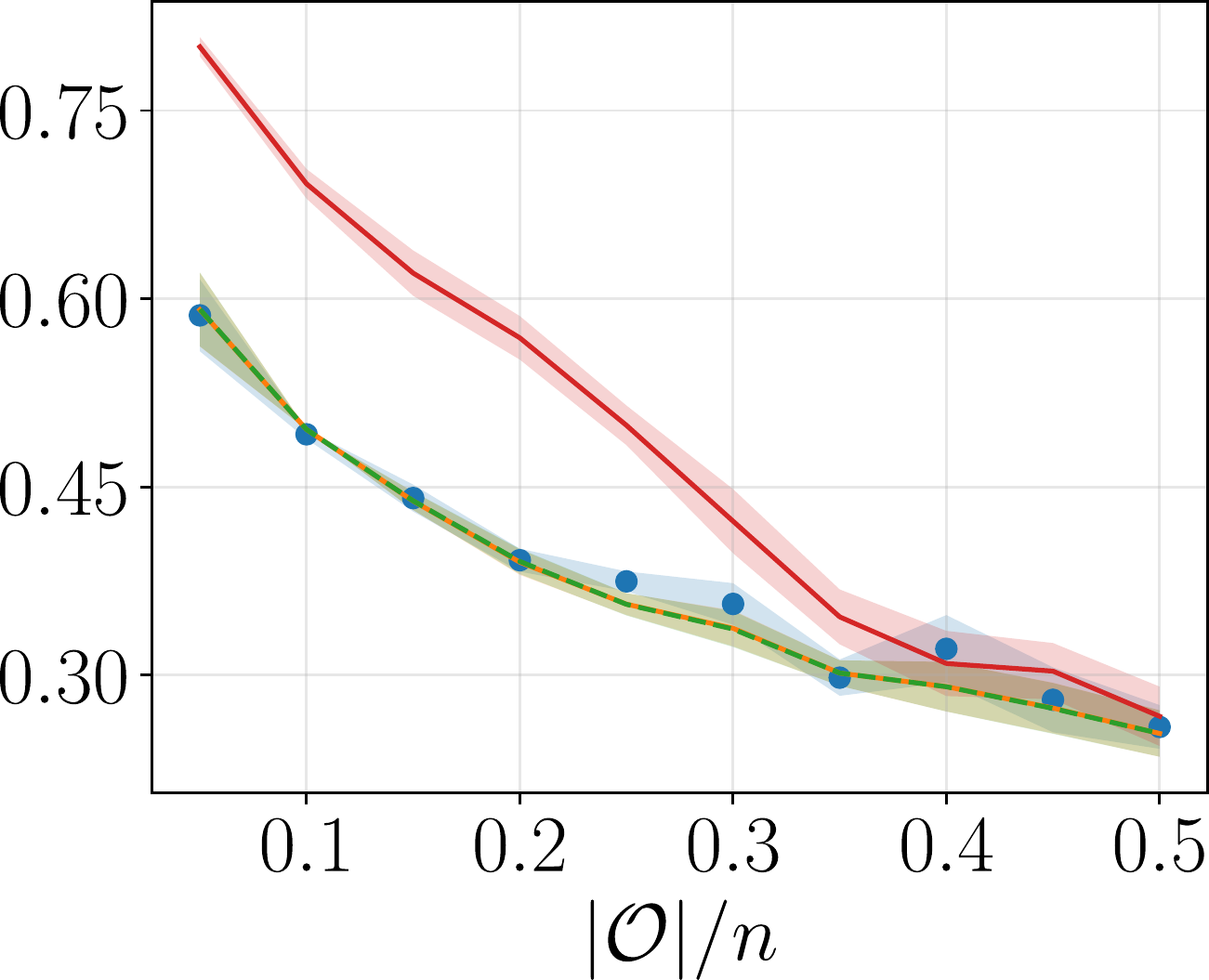}
		{\small (f)\; $\calO$: 7, $\calI$: all}
	\end{minipage}
	\begin{minipage}[t]{\ratio\linewidth}
		\centering
		\includegraphics[height=\ratiob\linewidth]{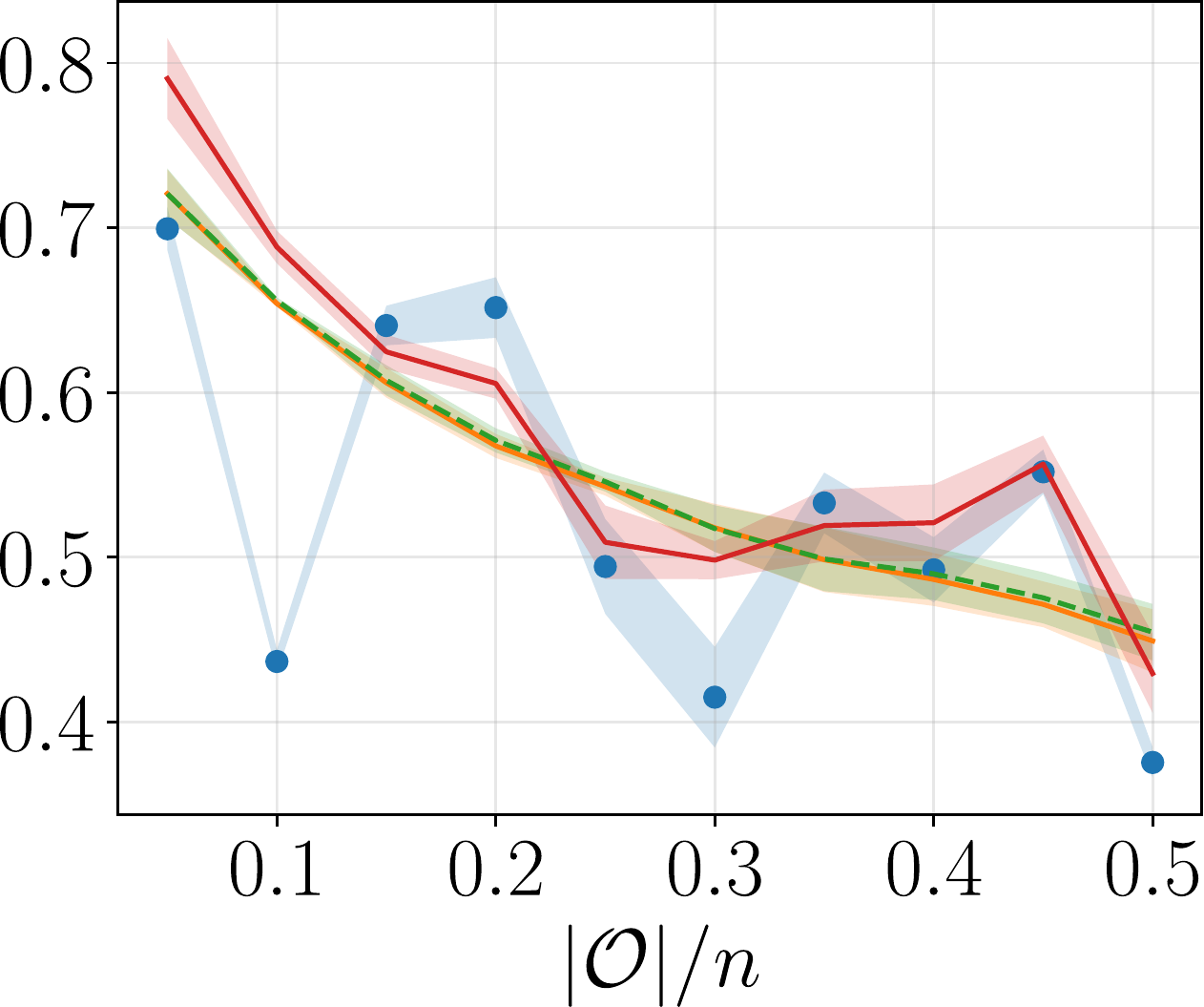}
		{\small (g)\; $\calO$: 8, $\calI$: all}
	\end{minipage}
	\begin{minipage}[t]{\ratio\linewidth}
		\centering
		\includegraphics[height=\ratiob\linewidth]{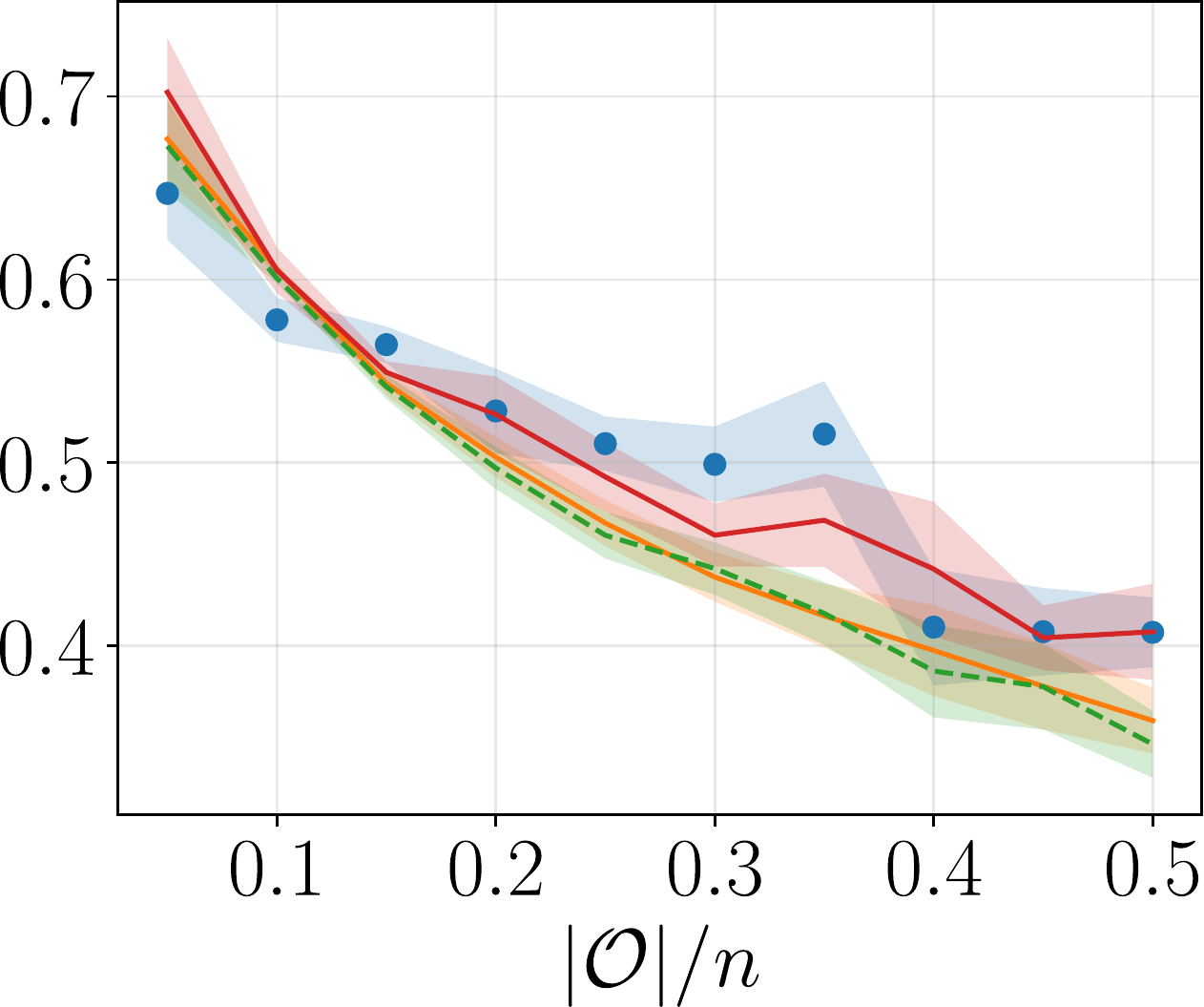}
		{\small (h)\; $\calO$: 9, $\calI$: all}
	\end{minipage}
	\caption{Anomaly detection with Digits data, measured with AUC over varying outlier proportion. A higher score means a better detection of the outliers. We specify which classes are chosen to be inliers ($\calI$) and outliers ($\calO$).}
	\label{fig:experiments_real_app}
\end{figure}

\end{document}